\documentclass[a4paper]{article}
\usepackage{amsmath,amssymb,fancybox,longtable,graphics,amscd,multirow,theorem}
\usepackage[dvipdfm]{graphicx} 
\usepackage{caption,setspace}
\newtheorem{thm}{\textsc{Theorem}}[section]

\newtheorem{lem}{\textsc{Lemma}}[section]
\newtheorem{cor}{\textsc{Corollary}}[section]

\def\problem{\textsc{Problem}\quad} 
\def\proof{\textsc{Proof. }}
\def\QED{$\Box$} \def\END{$\blacksquare$}
\def\Pic{\mathop{\mathrm{Pic}}\nolimits}
\def\discr{\mathop{\mathrm{discr}}\nolimits}
\def\rk{\mathop{\mathrm{rank}}\nolimits}

\def\sign{\mathop{\mathrm{sgn}}\nolimits}

\theorembodyfont{\rmfamily}
\newtheorem{remark}{\textsc{Remark}}

\theorembodyfont{\rmfamily}

\begin{document} 
\title{Lattice duality for coupling pairs admitting polytope
duality with trivial toric contribution}
\author{Makiko Mase}
\date{\empty}
\maketitle
\abstract{We study a lattice duality among families of $K3$ surfaces associated to coupling pairs that admit polytope duality with trivial toric contribution. } \\


\section{Introduction}
Weight systems appear in many interesting spots in algebraic geometry including singularity theory, where singularities have nice properties. 
We focus on a duality among weight systems called {\it coupling} introduced by Ebeling~\cite{EbelingCoupling}, which is for well-posed weight systems associated to simple $K3$ singularities classified by Yonemura~\cite{Yonemura}. 
The coupling duality is in particular admitted by a pair of singularities defined by weighted-homogeneous polynomials $f$ and $f'$ as a strange-duality for invertible polynomials introduced by Ebeling and Takahashi in ~\cite{EbelingTakahashi11}. 
It is also known that such polynomials $f$ and $f'$ in three variables can be projectivised as weighted-homogeneous polynomials $F$ and $F'$ as anticanonical divisor of the weighted projective spaces $\mathbb{P}_a$ and $\mathbb{P}_b$, where the pair $(a,\, b)$ is coupling among Yonemura's list. 
Since all the weighted projective spaces with weights being in Yonemura's list are Fano, we obtain subfamilies of $K3$ surfaces in the space once one finds a reflexive polytope as a subpolytope of the defining polytope of the space. 
In the author's recent work~\cite{Mase19-1}, an existence and duality of such reflexive polytopes are studied and it is concluded that almost all coupling pair extends to a polytope-duality.  
Once one obtains families of $K3$ surfaces which already admit several dualities, one may be interested in intrinsic properties of $K3$ surfaces. 
We are interested in lattice-duality originally studied by Dolgachev~\cite{DolgachevMirror}. 
It is concluded by the author~\cite{Mase15, Mase17} that a part of transpose-dual pairs associated to strange duality of bimodal singularities extends to lattice dual, and that some subfamilies of $K3$ surfaces that are double covering of the projective plane have lattice-dual property as is studied in~\cite{Mase19}. 
In this paper, focusing on polytope-dual pairs associated to coupling, one may pose the following problem. \\

\noindent
\problem
Determine whether or not the coupling pairs which admit polytope-duality extend to lattice duality of families $\mathcal{F}_{\Delta}$ and $\mathcal{F}_{\Delta'}$ in the sense that the relation
\[
(\Pic_{\Delta})^\perp_{U^{\oplus 3}\oplus E_8^{\oplus 2}} \simeq U\oplus \Pic_{\Delta'} 
\]
holds. 

We give an answer as the main theorem of the article which is presented here : \\

\bigskip
\noindent
{\bf Theorem \ref{MainThm}}\quad 
{\it If a coupling pair admits polytope-duality with trivial toric contribution, then, the families of $K3$ surfaces are lattice dual. 
Explicite Picard lattices of the families are given in Table~\ref{ListMainThm}. }
\\

In section $2$, we recall the Picard lattice and toric geometry. 
In section $3$, we give a proof of the main theorem. 
In the last and fourth section, we give a conclusion as the property of the Picard lattices of families that we have obtained.

\section{Preliminery}
A {\it lattice} is a finitely-generated $\mathbb{Z}$-module with a non-degenerate bilinear form. 
A {\it $K3$ surface} is a smooth compact complex connected $2$-dimensional algebraic variety with trivial canonical divisor and irregularity zero. 
It is known that the second cohomology group with $\mathbb{Z}$-coefficient of a $K3$ surface $S$ admits a structure of a unimodular lattice of signature $(3,19)$, thus by a classification of lattices, the lattice is in fact isometric to the {\it $K3$ lattice} $\Lambda_{K3}:=U^{\oplus 3}\oplus E_8^{\oplus 2}$, where $U$ is the hyperbolic lattice of rank $2$, and $E_8$ is the negative-definite, even unimodular lattice of rank $8$. 
By a standard exact sequence, one gets an inclusion map $c_1:H^{1}(S,\,\mathcal{O}_S^*)\to H^2(S,\,\mathbb{Z})$, which makes the Picard group $H^{1}(S,\,\mathcal{O}_S^*)$ to be a sublattice of $H^2(S,\,\mathbb{Z})$. 
We call the Picard group of $S$ with a lattice structure simply the {\it Picard lattice} of $S$. \\

We summerize toric geometry in \cite{BatyrevMirror} by also giving useful formulas extracted from \cite{FultonToric} and \cite{Oda78}. 

Let $M$ be a lattice of rank $n$, and $N:={\mathrm{Hom}}_{\mathbb{Z}}(M\,\mathbb{Z})$ be the dual lattice of $M$, with a natural pairing $\langle\, ,\,\rangle: N\times M\to\mathbb{Z}$ with its $\mathbb{R}$-extension denoted by $\langle\, ,\,\rangle_{\mathbb{R}}$. 
A convex hull of finite-number of points in $M\otimes\mathbb{R}$ is called a {\it polytope}, which admits the {\it polar dual} polytope $\Delta^*$ defined by 
\[
\Delta^*:=\left\{ y\in N\otimes\mathbb{R}\, |\, \langle y,x\rangle_{\mathbb{R}}\geq -1\quad \textnormal{for all } x\in\Delta\right\}. 
\]
A polytope $\Delta$ is {\it integral} if every vertex is in $M$. 
An integral polytope $\Delta$ which contains the only lattice point in its interior is {\it reflexive} if the polar dual polytope $\Delta^*$ is also an integral polytope. 

It is observed by \cite{BatyrevMirror} that an integral polytope $\Delta$ is reflexive if and only if the resulting projective toric variety $\mathbb{P}_{\!\Delta}$ is Fano, in other words, general hypersurfaces that are defined by global anticanonical sections of $\mathbb{P}_{\!\Delta}$ are birational to Calabi-Yau. 

We only treat with $3$-dimensional reflexive polytopes. 
We call a {\it anticanonical section} for hypersurfaces that are defined by global anticanonical sections of $\mathbb{P}_{\!\Delta}$ for short. 
In $3$-dimensional case, it is derived by a study of \cite{BatyrevMirror}, that moreover, singularities in $\mathbb{P}_{\!\Delta}$ and in general anticanonical sections $Z$ of $\mathbb{P}_{\!\Delta}$ can be simultaneously resolved by a toric resolution called a  MPCP-desingularisation, which we denote by $\widetilde{\mathbb{P}_{\!\Delta}}$ and $\tilde{Z}$. 
The natural restriction map 
\[
H^{1,1}(\widetilde{\mathbb{P}_{\!\Delta}},\,\mathbb{Z}) \to H^{1,1}(\tilde{Z},\,\mathbb{Z}) 
\]
is not necessarily surjective in general, and we denote by $L_{0}(\Delta)$ the rank of the cokernel of the map, which we call the {\it toric contribution}, which is known \cite{Kobayashi} to be given by the formula
\begin{equation}\label{ToricContribution}
L_0(\Delta) = \sum_{\Gamma}l(\Gamma)l(\Gamma^*),
\end{equation}
where the sum runs for all edges in $\Delta$. 

Here we recall from \cite{Bruzzo-Grassi} that generic anticanonical sections of the Fano $3$-fold $\mathbb{P}_{\!\Delta}$ admit isomemtric Picard lattices. 
Thus, we define the {\it Picard lattice of the family} $\mathcal{F}_{\!\Delta}$ of $K3$ surfaces in $\mathbb{P}_{\!\Delta}$ to be the Picard lattice of the minimal model of any generic anticanonical section of $\mathbb{P}_{\!\Delta}$, and denote it by $\Pic_{\Delta}$. 

For a reflexive polytope $\Delta$, one can associate a fan $\Sigma'$. 
By definition, lattice points of $\Delta^*$ are primitive vector of one-simplices of $\Sigma'$, and it is clear that the toric varieties $\mathbb{P}_{\!\Delta}$ and $\mathbb{P}_{\Sigma'}$ coincide. 
Any divisor $D$ of a generic hypersurface in $\mathbb{P}_{\!\Delta}$ is the closure of the torus orbit of a one-simplex $v$ in $\Sigma$, in particular, the divisors are called {\it toric divisors}. 
Let $F$ be the face in $\Delta$ that is the polar dual of $v$. 
Denote by $l(F)$ the number of lattice points in the interior of $F$. 
The self-intersection number of the divisor $D$ is given by the formula 
\begin{equation}\label{self-intersection}
D^2 = 2l(F)-2. 
\end{equation}
Denote by $\Delta^{(1)}$ the set of all edges in $\Delta$ and $l(\Gamma)$ be the number of lattice points in the interior of an edge $\Gamma\in\Delta^{(1)}$. 
The Picard number $\rho(\Delta)$ is given by 
\begin{equation}\label{PicardNumber}
\rho(\Delta) = \sum_{\Gamma\in\Delta^{(1)}}l(\Gamma) + \sum_{\textnormal{vertices of} \Delta}1+ L_{0}(\Delta)-3,
\end{equation}
Let $e_1,\, e_2,\, e_3$ be a standard basis for $\mathbb{R}^3$. 
Suppose that the fan $\Sigma$ possesses $l$ one-simplices. 
The toric divisors $D_1,\ldots,D_l$ admit the linear relations 
\begin{equation}\label{ToricDivisorsRelations}
\sum_{i=1}^l \langle v_i,\, e_j\rangle D_i= 0 \quad j=1,2,3. 
\end{equation}

It is easily seen that the polytope $\Delta$ is of trivial toric contribution if and only if the corresponding fan $\Sigma'$ is simplicial, that is, every triple of one-simplices form a $\mathbb{Z}$-basis of $\mathbb{R}^3$. 
Moreover, the restriction of linearly-independet toric divisors of $X=\widetilde{\mathbb{P}_{\!\Delta}}=\widetilde{\mathbb{P}_{\Sigma'}}$ to the anticanonical divisor of $X$ form a basis of the Picard lattice $\Pic_{\!\Delta'}$ of the family $\mathcal{F}_{\!\Delta'}$ if $\Delta'\simeq\Delta^*$. 

Denote by $M_{(a_0, a_1, a_2, a_3)}$ the lattice consisting of quadruple of integers $(i, j, k, l)$ satisfying an equation $a_0i+a_1j+a_2k+a_3l=0$ for a weight system $(a_0, a_1, a_2, a_3; d)$. 
There is a one-to-one correspondence between elements in $M_{(a_0, a_1, a_2, a_3)}$ and (rational) monomials of degree $d$ by 
\[
(i, j, k, l)\in M_{(a_0, a_1, a_2, a_3)} \leftrightarrow W^{i+1}X^{j+1}Y^{k+1}Z^{l+1},  
\]
where $(W,X,Y,Z)$ is a coordinate system of the weighted projective space of weight $(a_0, a_1, a_2, a_3)$. 

\bigskip
We denote by $L^*$,\, $A_L$,\, $\discr{L}$, \, $l(A_L)$, \, $\sign{L}$,\, $q_L$, and $\rk{L}$ the {\it dual lattice} $L^*:={\mathrm Hom}_{\mathbb{Z}}(L,\,\mathbb{Z})$, the {\it discriminang group} $L\slash L^*$, the {\it discriminant}, the {\it minimal number of generators of $A_L$}, the {\it signature}, the {\it discriminant form}, and the {\it rank} of a lattice $L$. 
It is a standard arithmetic property that if $\rk{L}$ is strictly larger than $5$, then, there eists an element representing $0$, and if $\rk{L}$ is strictly larger than $12$, then, the hyperbolic lattice $U$ is a sublattice of $L$. 
We also recall standard properties of lattices from \cite{Nikulin80} and \cite{Nishiyama96}. 
A sublattice $S$ of a lattice $\Lambda$ is called {\it primitive} if the quotient lattice $\Lambda\slash S$ is torsion-free. 

\begin{cor}{\rm (Corollary 1.6.2~\cite{Nikulin80})}\label{NikulinOrthogonal}
Let $S$ and $T$ be primitive sublattices of the $K3$ lattice $\Lambda_{K3}$. 
The lattices $S$ and $T$ are orthogonal in $\Lambda_{K3}$ if and only if $q_S\simeq-q_T$ holds. \END
\end{cor}
\begin{cor}{\rm (Corollary 1.12.3~\cite{Nikulin80})}\label{NikulinPrimitive}
Let $S$ be a sublattice with signature $(t_+,\, t_-)$ of an even unimodular lattice $\Lambda$ with signature $(l_+,\, l_-)$. 
The lattice $S$ is a primitive sublattice of $\Lambda$ if and only if the following three conditions are satisfied. \\
$(1)$\, $l_+ - l_- \equiv 0 \mod 8$, \\
$(2)$\, $l_- - t_-\geq 0$ and $l_+ - t_+\geq 0$, and \\
$(3)$\, $\rk{\Lambda}-\rk{S} >l(A_S)$. 
\END
\end{cor}
\begin{remark}
Note that the $K3$ lattice $\Lambda_{K3}$ is an even unimodular lattice of signature $(l_+,\, l_-)=(3,\, 19)$. 
Thus, $l_+ - l_- = 3-19=-16\equiv0\mod 8$, and in order to show a lattice $S$ to be a primitive sublattice of $\Lambda_{K3}$, it suffices to verify the second and third conditions of Corollary~\ref{NikulinPrimitive}. 
\end{remark}
\begin{lem}{\rm (Lemma 4.3~\cite{Nishiyama96})}\label{Nishiyama}
There exist primitive embeddings of $A_1$ and $A_2$ into $E_8$ with orthogonal complements being $E_7$ and $E_6$, respectively. 
We follow the notation of lattices in Bourbaki~\cite{Bourbaki68}. 
\END
\end{lem}

\section{Main Results}
\begin{lem}
The polytope-dual pairs among Nos. 11-14, Nos. 15-18, Nos. 35--37, Nos. 38 and 40, Nos. 41--43, Nos. 48-49 are respectively isomorphic to the following polytopes in Table~\ref{ListMainThm}. 
\end{lem}
\proof
The assertion follows from the proof of \cite{Mase19-1}. 
\QED

\begin{lem}
If a coupling pair is in Talbe~\ref{ListMainThm}, the toric contribution is trivial. 
\end{lem}
\proof
The assertion follows by case-by-case computation using formula (\ref{ToricContribution}) for all polytopes obtained in \cite{Mase19-1}. 
\QED

\begin{thm}\label{MainThm}
If a coupling pair admits polytope-duality with trivial toric contribution, then, the families of $K3$ surfaces are lattice dual. 
Explicite Picard lattices of the families are given in Table~\ref{ListMainThm}. 
\begingroup
\setlength{\arrayrulewidth}{.1pt}
\tiny
\begin{longtable}[htp]{p{1.8mm}|
p{3mm}
p{30mm}@{\hspace{7mm}}
p{10mm}@{\hspace{7mm}}
p{10mm}@{\hspace{10mm}}
p{35mm}}
{\rm No.} & & \hspace{10mm}$\Delta'$ & $\Pic(\Delta')$, $(\rk, \lvert\discr\rvert)$, \textnormal{weight system} & $\Pic(\Delta)$, $(\rk, \lvert\discr\rvert)$, \textnormal{weight system} &\hspace{15mm} $\Delta$ \\
\hline
\shortstack{\\$11$,\\$12$,\\$13$,\\$14$}
& \raisebox{1.8mm}{\shortstack{Case 1 \\ \dotfill \\ Case 2}}
&  \raisebox{1mm}{\shortstack{$Z^2,\, W^{30},\, W^6X^6,$ $X^5Y,\, Y^3$ \\ \dotfill \\ $Z^2,\, W^{30},\, W^2X^7,$ $X^5Y,\, Y^3$} }
& \raisebox{1mm}{\shortstack{$U\oplus E_7$, \\ $(9,2)$ \\ $1,4,10,15;30$}}
& \raisebox{1mm}{\shortstack{$U\oplus A_1\oplus E_8$, \\  $(11,2)$ \\ $1,6,8,15;30$}}
& \raisebox{1mm}{\shortstack{$Z^2,\, W^{30},\, X^5,$ $XY^3,\, W^6Y^3$ \\ \dotfill \\ $Z^2,\, W^{30},\, X^5,$ $XY^3,\, W^{14}Y^2$}} \\
\hline
\hline
\shortstack{\\$15$,\\$16$,\\$17$,\\$18$}
&
& \raisebox{4mm}{\shortstack{\\$Z^2,\, W^{24},\, W^6X^6,$ $X^4Z,\, Y^3$ }}
&  \raisebox{1mm}{\shortstack{\\$U\oplus E_6$, \\ $(8,3)$ \\ $ 1,3,8,12;24$}}
&  \raisebox{1mm}{\shortstack{\\$U\oplus A_2\oplus E_8$, \\ $(12,3)$ \\ $1,6,8,9;24$}}
&  \raisebox{4mm}{\shortstack{\\$W^6Z^2,\, W^{24},\, X^4,$ $XZ^2,\, Y^3$ }}
\\
\hline
\hline
\raisebox{10mm}{$19.$}
& \raisebox{0.8mm}{\shortstack{Case 1 \\ \\ \dotfill \\ \\ \\ \\ Case 2 \\ \\ \dotfill \\ \\ \\ \\  Case 3 }}
& \shortstack{\\$Z^2,\, W^{22},\, W^2X^5,$\\ $X^4Y,\, XY^3,\, W^{10}Y^2$ \\ \dotfill 
\\ $Z^2,\, W^{22},\, W^6X^4,$\\ $X^4Y,\, XY^3,\, W^4Y^3$ \\ \dotfill 
\\ $Z^2,\, W^{22},\, W^6X^4,$\\ $X^4Y,\, XY^3,\, W^{10}Y^2$} 
& \raisebox{8mm}{\shortstack{\\$ U\oplus A_1\oplus E_7$, \\ $(10, 4)$ \\ $1,4,6,11;22$}}
& \raisebox{8mm}{\shortstack{\\$ U\oplus A_1\oplus E_7$, \\ $(10, 4)$ \\ $1,4,6,11;22$}}
& \shortstack{\\$Z^2,\, W^{22},\, W^2X^5,$\\ $X^4Y,\, XY^3,\, W^{10}Y^2$  \\ \dotfill 
\\ $Z^2,\, W^{22},\, W^6X^4,$\\ $X^4Y,\, XY^3,\, W^4Y^3 $\\ \dotfill 
\\ $Z^2,\, W^{22},\, W^2X^5,$\\ $X^4Y,\, XY^3,\, W^4Y^3$}
\\
\hline
\hline
\raisebox{14mm}{$26.$}
& \raisebox{0.8mm}{\shortstack{Case 1 \\ \\ \dotfill \\ \\ \\ \\ Case 2 \\ \\ \dotfill \\ \\ \\ \\  Case 3 \\ \\ \dotfill \\ \\ \\ \\ Case 4}}
& \shortstack{\\ $W^3Z^2,\, Y^2Z,\, XZ^2,\, W^{13},$\\ $W^4X^3,\, X^3Y,\,WY^3$ \\ \dotfill \\ $W^8Z,\, Y^2Z,\, XZ^2,\,W^{13},$\\ $WX^4,\, X^3Y,\,WY^3$ \\ \dotfill \\ $W^3Z^2,\, Y^2Z,\, XZ^2,\,W^{13},$\\ $WX^4,\, X^3Y,$ $W^9Y$  \\ \dotfill \\ $W^3Z^2,\, Y^2Z,\, XZ^2,\,W^{13},$\\ $WX^4,\, X^3Y,\,WY^3$}
& \raisebox{14mm}{\shortstack{$U\oplus \tilde{L'}$,\\ $(10,13)$ \\ $1,3,4,5;13$}}
& \raisebox{14mm}{\shortstack{$U\oplus \tilde{L}$,\\ $(10,13)$ \\ $1,3,4,5;13$}}
& \shortstack{\\ $W^3Z^2,\, Y^2Z,\, XZ^2,\,W^{13},$\\ $W^4X^3,\, X^3Y,$ $W^9Y$ \\ \dotfill \\ $W^8Z,\, Y^2Z,\, XZ^2,\,W^{13},$\\ $WX^4,\, X^3Y,$ $W^9Y$ \\ \dotfill \\  $W^8Z,\, Y^2Z,\, XZ^2,\,W^{13},$\\ $W^4X^3,\, X^3Y,\,WY^3$ \\ \dotfill \\ $W^8Z,\, Y^2Z,\, XZ^2,\,W^{13},$\\ $W^4X^3,\, X^3Y,$ $W^9Y$}
\\
\hline
\hline
\shortstack{\\$35$,\\$36$,\\$37$}
&
&  \raisebox{3mm}{\shortstack{\\$Z^2,\, W^{12},\, X^{12},\, Y^3$}}
& \shortstack{$U$,\\ $(2,1)$ \\ $1,1,4,6;12$}
& \shortstack{\\$U\oplus U\oplus E_8^{\oplus2}$,\\ $(18,1)$ \\ $3,5,11,14;33$}
&   \raisebox{3mm}{\shortstack{\\$XZ^2,\, W^{11},\, WX^6,\, Y^3$}}\\
\hline
\hline
\raisebox{1.4mm}{\shortstack{\\$38$,\\$40$}}
& \raisebox{0.8mm}{\shortstack{Case 1 \\ \dotfill \\ Case 2}}
& \shortstack{\\$Z^2,\, W^{10},\, X^{10},\,XY^3,\, WY^3$ \\ \dotfill \\ $Z^2,\, W^{10},\, X^{10},\,XY^3,\, W^4Y^2$}
& \shortstack{$U\oplus E_7\oplus E_8$,\\ $(17,2)$ \\ $1,1,3,5;10$}
& \shortstack{$U\oplus A_1$,\\ $(3,2)$ \\ $3,4,10,13;30$}
& \shortstack{\\$XZ^2,\, W^{10},\, W^6X^3,\,X^5Y,\, Y^3$ \\ \dotfill \\ $XZ^2,\, W^{10},\, W^2X^6,\,X^5Y,\, Y^3 $}
\\
\hline
\hline
\shortstack{\\$41$,\\$42$,\\$43$}
&
& \raisebox{3mm}{\shortstack{\\Y$^3,\, WZ^2,\, W^9,\,X^9,\, XZ^2$}}
& \shortstack{\\$U\oplus A_2$, \\ $(4,3)$ \\ $1,1,3,4;9$}
& \shortstack{\\$U\oplus \oplus E_6\oplus E_8$, \\ $(16,3)$ \\ $3,4,11,18;36$}
& \raisebox{3mm}{\shortstack{$WY^3,\, Z^2,\, W^6Z,\,W^4X^6,\, X^9$}}
\\
\hline
\hline
\raisebox{3mm}{$46.$}
&
& \raisebox{3mm}{\shortstack{\\$Y^3Z,\, WZ^2,\, X^3Z,\, W^5,\, X^5,\, Y^5$}}
& \shortstack{\\$\left(\begin{smallmatrix} 2 & 1 \\ 1 & -2 \end{smallmatrix}\right)$, \\ $(2,5)$ \\ $1,1,1,2;5$}
& \shortstack{$U\oplus \tilde{L}$, \\ $(18,5)$ \\ $4,5,7,9;25$}
& \raisebox{3mm}{\shortstack{\\$YZ^2,\, W^4Z,\, W^5X,\, WY^3,\, X^5$}}
\\
\hline
\hline
\raisebox{2mm}{\shortstack{$48$, \\$49$}}
&
&  \raisebox{3mm}{\shortstack{\\$Z^2,\, W^6,\, X^6,\, Y^6$}}
& \shortstack{\\$\langle 2\rangle$, \\ $(1,2)$ \\ $1,1,1,3;6$}
& \shortstack{\\$U\oplus \langle -2\rangle\oplus E_8^{\oplus 2}$, \\ $(19,2)$ \\ $5,6,8,11;30$}
& \raisebox{3mm}{\shortstack{\\$YZ^2,\, W^6,\, X^5,\, XY^3$}}\\
\hline
\hline
\raisebox{2mm}{$50.$}
&
& \raisebox{3mm}{\shortstack{\\$Z^4,\, W^4,\, X^4,\, Y^4$}}
& \shortstack{\\$\langle 4\rangle$, \\ $(1,4)$ \\ $1,1,1,1;4$}
& \shortstack{\\$U\oplus \langle -4\rangle\oplus E_8^{\oplus 2}$, \\ $(19,4)$ \\ $7,8,9,12;36$}
& \raisebox{3mm}{\shortstack{\\$Z^3,\, W^4X,\, X^3Z,\, Y^4$}}\\
\hline
\hline
\hline
\caption{Lattice duality associated to coupling pairs}\label{ListMainThm}
\end{longtable}
\endgroup
\end{thm}

\begin{remark}
We present the following data in Table~\ref{ListMainThm}. 
The number(s) in the first column are given in \cite{Mase19-1}. 
The second and fifth columns are vertices of polytopes of $\Delta'$ and $\Delta$ obtained by \cite{Mase19-1}, and the sets in the same line are polytope-dual. 
In the third and fourth columns are the Picard lattice of the family $\mathcal{F}_{\!\Delta'}$, resp. $\mathcal{F}_{\!\Delta}$, the pair of the rank and the  signature of lattices, and the weight systems that are coupling. 
The lattices $\tilde{L}$ and $\tilde{L'}$ are explained in the proof. 
Other lattices follow notation of \cite{Bourbaki68}. 
\end{remark}

\noindent
\proof
Take reflexive polytopes $\Delta$ and $\Delta'$ as in \ref{ListMainThm}. 
We explicitely calculate the Picard lattices of the families $\mathcal{F}_{\!\Delta}$ and $\mathcal{F}_{\!\Delta'}$. 
Denote by $\Sigma$, respectively $\Sigma'$ the fan associated to polytope $\Delta'$, resoectively $\Delta$. 
Since the relation $\Delta^*\simeq\Delta'$ holds, lattice points of $\Delta'$, respectively of $\Delta$ are none other than primitive vectors of one-simplices of $\Sigma$, respectively $\Sigma'$. 
\subsection{Nos. 11, 12, 13, and 14}
\underline{Case 1.}\quad 
Set one-simplices of $\Sigma$ in terms of a basis of $M_{(1,4,10,15)}\otimes\mathbb{R}$ 
\[
(-1,4,0,-1),\, (-1,-1,2,-1),\, (-1,-1,-1,1) : 
\]
\[
\begin{array}{lll}
v_1 = (1, 0, 0), & v_2 = (0, 1, 0), & v_3 = (0, 0, 1), \\
v_4 = (0, -2, -3), & v_5 = (-6, -8, -15), & v_6 = (-3, -4, -7), \\
v_7 = (0, -1, -1), & v_8 = (-2, -2, -5), & v_9 = (-4, -5, -10), \\
v_{10} = (-5, -7, -13), & v_{11} = (-4, -6, -11), & v_{12} = (-3, -5, -9), \\
v_{13} = (-2, -4, -7), & v_{14} = (-1, -3, -5), 
\end{array}
\]
and let $\tilde{D_i}$ be the toric divisor determined by the lattice point $v_i$ for $i=1,\ldots,14$, and $D_i:=\tilde{D_i}|_{{-}K_X}$ with $X:=\widetilde{\mathbb{P}_{\Sigma}}$. 
It can be easily seen by formulas (\ref{PicardNumber}) and (\ref{self-intersection}) that 
\[
\rho_{\Delta} = 14-3 = 11, \quad
D_1^2 = 0, \, D_2^2 = 2, \, D_3^2 = 8, \, D_4^2 = \cdots = D_{14}^2 = -2. 
\]
Let $L$ be a lattice generated by divisors $\{D_i \}_{i=1}^{14}$. 
By solving the equation (\ref{ToricDivisorsRelations}), one sees that $\{$ $D_1$, $D_2$, $D_8$, $D_4$, $D_7$, $D_{14}$, $D_{13}$, $D_{12}$, $D_{11}$, $D_{10}$, $D_5$ $\}$ form a basis for $L$. 
By taking a new basis
\[
\left\{
\begin{array}{l}
D_1,\, D_1+D_4,\, -D_1+D_7,\, D_8,\\
2D_1-D_2+2D_4+D_7-D_8+D_{14}, \\
D_{13},\, -D_1+D_{14},\, D_{12},\, D_{11},\, D_{10},\, D_5
\end{array}
\right\}, 
\]
one sees that the lattice $L$ is isometric to $U\oplus A_1\oplus E_8$, which is hyperbolic and a primitive sublattice of the $K3$ lattice. 
Therefore, $\Pic_{\!\Delta}\simeq U\oplus A_1\oplus E_8$. 

Set one-simplices of $\Sigma'$ in terms of a basis of $M_{(1,6,8,15)}\otimes\mathbb{R}$
\[
(-6,1,0,0),\, (-8,0,1,0),\, (-15,0,0,1) :
\]
\[
\begin{array}{lll}
m_1 = (4, -1, -1), & m_2 = (0, 2, -1), & m_3 = (-1, -1, 1), \\
m_4 = (-1, 2, -1), & m_5 = (-1, -1, -1), & m_6 = (3, -1, -1), \\
m_7 = (2, -1, -1), & m_8 = (1, -1, -1), & m_9 = (0, -1, -1), \\
m_{10} = (-1, 1, -1), & m_{11} = (-1, 0, -1), & m_{12} = (-1, -1, 0),
\end{array}
\]
and let $\tilde{D_i'}$ be the toric divisor determined by the lattice point $m_i$ for $i=1,\ldots,12$, and $D_i':=\tilde{D_i'}|_{{-}K_X}$ with $X:=\widetilde{\mathbb{P}_{\Sigma'}}$. 
It can be easily seen by formulas (\ref{PicardNumber}) and (\ref{self-intersection}) that 
\[
\rho_{\Delta'} = 12-3 = 9, \quad
D_1'^2 = 0, \, D_2'^2 = 2, \, D_3'^2 = 10, \, D_4'^2 = \cdots = D_{12}'^2 =-2. 
\]
Let $L$ be a lattice generated by divisors $\{D_i' \}_{i=1}^{12}$. 
By solving the equation (\ref{ToricDivisorsRelations}), one sees that $\{ D_1', D_2', D_4', D_{10}', D_{11}',  D_5', D_6', D_7', D_8' \}$ form a basis for $L$. 
By taking a new basis
\[
\left\{
\begin{array}{l}
D_1', \, D_1'+D_6', \, D_4',\, D_{10}',\, D_{11}',\, D_5', \\
-3D_1'+D_2'+D_4'-2 D_6'-D_7'-D_8',\, -D_1'+D_7'+D_8',\, D_1'-D_7'
\end{array}
\right\}, 
\]
one sees that the lattice $L'$ is isometric to $U\oplus E_7$, which is hyperbolic and a primitive sublattice of the $K3$ lattice. 
Therefore, $\Pic_{\!\Delta'}\simeq U\oplus E_7$. \\

\noindent
\underline{Case 2}\quad 
Set one-simplices of $\Sigma$ in terms of a basis of $M_{(1,4,10,15)}\otimes\mathbb{R}$ 
\[
(-4,1,0,0),\, (-10,0,1,0),\, (-15,0,0,1) :
\]
\[
\begin{array}{lll}
v_1 = (-1, -1,  1), & v_2 = (-1, -1,  -1), & v_3 =  (6,  -1,  -1), \\
v_4 =  (4,  0,  -1), & v_5 =  (-1,  2,  -1), & v_6 =  (-1, -1,  0), \\
v_7 =  (0,  -1,  -1), & v_8 =  (1,  -1,  -1), & v_9 =  (2,  -1,  -1), \\
v_{10} = (3,  -1,  -1), & v_{11} = (4,  -1,  -1), & v_{12} = (5,  -1,  -1), \\
v_{13} = (-1,  0,  -1), & v_{14} = (-1,  1,  -1), 
\end{array}
\]
and let $\tilde{D_i}$ be the toric divisor determined by the lattice point $v_i$ for $i=1,\ldots,14$, and $D_i:=\tilde{D_i}|_{{-}K_X}$ with $X:=\widetilde{\mathbb{P}_{\Sigma}}$. 
It can be easily seen by formulas (\ref{PicardNumber}) and (\ref{self-intersection}) that 
\[
\begin{array}{llll}
\rho_{\Delta} = 14-3 = 11, & 
D_1^2 = 8, \, D_2^2 = D_3^2 = D_4^2 = -2, \, D_5^2 = 2, \\
 &  D_6^2 = \cdots = D_{14}^2 =-2. 
\end{array}
\]
Let $L$ be a lattice generated by divisors $\{D_i \}_{i=1}^{14}$. 
By solving the equation (\ref{ToricDivisorsRelations}), one sees that $\{ D_4, D_3, D_5, D_{14}, D_{13}, D_ 2, D_{12}, D_{11}, D_{10}, D_9, D_8 \}$ form a basis for $L$. 
By taking a new basis
\[
\left\{
\begin{array}{l}
-D_3+D_5-D_{12},\, D_3+D_4,\, D_4,\, D_{13},\\
-D_2-D_{13},\, D_2-D_3-D_4+D_{13}+D_{14},\\
2D_3+D_4-D_5+D_{11}+2D_{12}-D_{13}-D_{14},\\ 
D_3+D_4+D_{11},\, D_{10},\, D_9,\, D_8
\end{array}
\right\}, 
\]
one sees that the lattice $L$ is isometric to $U\oplus A_1\oplus E_8$, which is hyperbolic and a primitive sublattice of the $K3$ lattice. 
Therefore, $\Pic_{\!\Delta}\simeq U\oplus A_1\oplus E_8$. \\

Set one-simplices of $\Sigma'$ in terms of a basis of $M_{(1,6,8,15)}\otimes\mathbb{R}$ 
\[
(-6,1,0,0),\, (-8,0,1,0),\, (-15,0,0,1) :
\]
\[
\begin{array}{lll}
m_1=( -1,  -1, 1), & m_2=( -1, -1,  -1), & m_3=( 4,  -1,  -1), \\
m_4=( 0,  2,  -1), & m_5=( -1,  1,  -1), & m_6=( -1,  -1,  0), \\
m_7=( 0,  -1,  -1), & m_8=( 1, -1,  -1), & m_9=( 2, -1, -1), \\
m_{10}=( 3,  -1,  -1), & m_{11}=( -1,  0,  -1), & m_{12}=( -1,  0,  0), 
\end{array}
\]
and let $\tilde{D_i'}$ be the toric divisor determined by the lattice point $m_i$ for $i=1,\ldots,12$, and $D_i':=\tilde{D_i'}|_{{-}K_X}$ with $X:=\widetilde{\mathbb{P}_{\Sigma'}}$.
It can be easily seen by formulas (\ref{PicardNumber}) and (\ref{self-intersection}) that 
\[
\begin{array}{llll}
\rho_{\Delta'} = 12-3 = 9, & 
D_1'^2 = 10, \, D_2'^2 = -2, \, D_3'^2 = 0, \, D_4'^2 = 4, \\
 & D_5'^2 = \cdots = D_{12}'^2 =-2. 
\end{array}
\]
Let $L'$ be a lattice generated by divisors $\{D_i' \}_{i=1}^{12}$. 
By solving the equation (\ref{ToricDivisorsRelations}), one sees that $\{ D_3', D_{10}', D_8', D_6', D_7', D_2', D_{11}', D_5', D_{12}' \}$ form a basis for $L'$, with respect to which the intersection matrix of $L'$ is $U\oplus E_7$, which is hyperbolic and a primitive sublattice of the $K3$ lattice. 
Therefore, $\Pic_{\!\Delta'}\simeq U\oplus E_7$. 

It is well-known that lattices $U\oplus A_1\oplus E_8$ and $U\oplus E_7$ are primitive sublattices of the $K3$ lattice $\Lambda_{K3}$. 
Moreover, by Lemma~\ref{Nishiyama}, the relation 
$(\Pic_{\!\Delta})^\perp_{\Lambda_{K3}}\simeq U^\perp_{U^{\oplus3}}\oplus(A_1)^\perp_{E_8}\oplus (E_8)^\perp_{E_8}=U^{\oplus2}\oplus E_7\simeq U\oplus\Pic_{\!\Delta'}$ holds. 

\subsection{Nos. 15, 16, 17, and 18}
Set one-simplices of $\Sigma$ in terms of a basis of $M_{(1,6,8,9)}\otimes\mathbb{R}$ 
\[
(-6,1,0,0),\, (-8,0,1,0),\,(-9,0,0,1) :
\]
\[
\begin{array}{lll}
v_1 = ( -1, 2,  -1), & v_2 = ( -1,  -1,  -1), & v_3 = ( 5,  -1,  -1), \\
v_4 = ( 3,  -1,  0), & v_5 = ( -1,  -1,  1), & v_6 = ( -1,  1,  -1), \\
v_7 = ( -1,  0,  -1), & v_8 = ( 0,  -1,  -1), & v_9 = ( 1,  -1,  -1), \\
v_{10} = ( 2,  -1,  -1), & v_{11} = ( 3,  -1, -1), & v_{12} = ( 4,  -1, -1), \\
v_{13} = ( -1,  -1,  0), & v_{14} = ( 3,  0,  -1), & v_{15} = ( 1,  1,  -1),
\end{array}
\]
and let $\tilde{D_i}$ be the toric divisor determined by the lattice point $v_i$ for $i=1,\ldots,15$, and $D_i:=\tilde{D_i}|_{{-}K_X}$ with $X:=\widetilde{\mathbb{P}_{\Sigma}}$. 
It can be easily seen by formulas (\ref{PicardNumber}) and (\ref{self-intersection}) that 
\[
\begin{array}{llll}
\rho_{\Delta} = 15-3 = 12, & 
D_1^2 = 2, \, D_2^2 = D_3^2 = -2, \\ 
& D_4^2 = 0, \, D_5^2 = 4, \, D_6^2 = \cdots = D_{15}^2 = -2. 
\end{array}
\]
Let $L$ be a lattice generated by divisors $\{D_i \}_{i=1}^{15}$. 
By solving the equation (\ref{ToricDivisorsRelations}), one sees that $\{$ $D_4$, $D_5$, $D_{13}$, $D_2$,$D_7$, $D_3$, $D_{14}$, $D_{15}$, $D_{12}$, $D_{11}$, $D_{10}$, $D_9$ $\}$ form a basis for $L$. 
By taking a new basis
\[
\left\{
\begin{array}{l}
D_4,\, D_3+D_4, \, D_4-D_{14}-D_{15},\, -D_4+D_{14}, \, -D_7, \\
-D_{13},\, -D_4+D_{11}+D_{12},\\
3D_3+2D_4-D_5+D_{10}+2D_{11}+2D_{12}+2D_{14}+D_{15},\\
-D_2,\, -D_9-D_{10}-D_{11},\, D_{10},\, D_9
\end{array}
\right\}, 
\]
one sees that the lattice $L$ is isometric to $U\oplus A_2\oplus E_8$, which is hyperbolic and a primitive sublattice of the $K3$ lattice. 
Therefore, $\Pic_{\!\Delta}\simeq U\oplus A_2\oplus E_8$. 

Set one-simplices of $\Sigma'$ in terms of a basis of $M_{(1,3,8,12)}\otimes\mathbb{R}$
\[
(-3,1,0,0),\, (-8,0,1,0),\, (-12,0,0,1): 
\]
\[
\begin{array}{lll}
m_1 = (-1,  2, -1), & m_2 = ( -1, -1, -1), & m_3 = ( 3, -1, -1), \\
m_4 = ( 0, -1, 1), & m_5 = ( -1, -1, 1), & m_6 = ( -1, 1, -1), \\
m_7 = ( -1, 0, -1), & m_8 = ( 0, -1, 1), & m_9 = ( 1, -1, -1), \\
m_{10} = ( 2, -1, -1), & m_{11} = ( -1, -1, 0), 
\end{array}
\]
and let $\tilde{D_i'}$ be the toric divisor determined by the lattice point $m_i$ for $i=1,\ldots,11$, and $D_i':=\tilde{D_i'}|_{{-}K_X}$ with $X:=\widetilde{\mathbb{P}_{\Sigma'}}$. 
It can be easily seen by formulas (\ref{PicardNumber}) and (\ref{self-intersection}) that 
\[
\begin{array}{llll}
\rho_{\Delta'} = 11-3 = 8, & 
D_1'^2 = 4, \, D_2'^2 = -2, \, D_3'^2 = 0, \\ 
 & D_4'^2 = 6, \, D_5'^2 = -2, \, D_6'^2 = \cdots = D_{11}'^2 =-2. 
\end{array}
\]
Let $L'$ be a lattice generated by divisors $\{D_i' \}_{i=1}^{11}$. 
By solving the equation (\ref{ToricDivisorsRelations}), one sees that $\{$ $D_4'$, $D_3'$, $D_{10}'$, $D_9'$, $D_5'$, $D_{11}'$, $D_2'$, $D_7'$ $\}$ form a basis for $L'$. 
By taking a new basis
\[
\left\{
\begin{array}{l}
D_4'-D_3',\, 2D_3'-D_5'+D_9'+2D_{10}',\, -D_3'+D_5',\\ 
3D_3'-2D_5'+D_9'+3D_{10}'-D_{11}',\\
D_9',\, -D3'+D_4'+D_{11}',\\ 
D_2'-3D_3'+2D_5'-2D_9'-3D_{10}'+D_{11}',\, D_7'
\end{array}
\right\}, 
\]
one sees that the lattice $L'$ is isometric to $U\oplus E_6$, which is hyperbolic and a primitive sublattice of the $K3$ lattice. 
Therefore, $\Pic_{\!\Delta'}\simeq U\oplus E_6$. 

It is well-known that lattices $U\oplus A_2\oplus E_8$ and $U\oplus E_6$ are primitive sublattices of the $K3$ lattice $\Lambda_{K3}$. 
Moreover, by Lemma~\ref{Nishiyama}, the relation 
$(\Pic_{\!\Delta})^\perp_{\Lambda_{K3}}\simeq U^\perp_{U^{\oplus3}}\oplus(A_2)^\perp_{E_8}\oplus (E_8)^\perp_{E_8}=U^{\oplus2}\oplus E_6\simeq U\oplus\Pic_{\!\Delta'}$ holds. 

\subsection{No. 19}
In all cases, we set one-simplices in fans in terms of a basis of $M_{(1,4,6,11)}\otimes\mathbb{R}$
\[
(-4,1,0,0), (-6,0,1,0), (-11,0,0,1). 
\]

\noindent
\underline{Case 1}\quad
We have $\Delta\simeq\Delta'$. 

Set one-simplices of $\Sigma'$ as follows: 
\[
\begin{array}{lll}
m_1=( -1,  -1, 1), & m_2=( -1,  -1, -1), & m_3=( 4,  -1,  -1), \\
m_4=( 3,  0,  -1), & m_5=( 0,  2,  -1), & m_6=( -1, 1,  -1), \\
m_7=( -1, -1,  0), & m_8=( 0,  -1,  -1), & m_9=( 1,  -1,  -1), \\
m_{10}=( 2,  -1,  -1), & m_{11}=( 3,  -1,  -1), & m_{12}=( -1,  0,  0), \\
m_{13}=( -1,  0,  -1), 
\end{array}
\]
and let $\tilde{D_i'}$ be the toric divisor determined by the lattice point $m_i$ for $i=1,\ldots,13$, and $D_i':=\tilde{D_i'}|_{{-}K_X}$ with $X:=\widetilde{\mathbb{P}_{\Sigma'}}$. 
It can be easily seen by formulas (\ref{PicardNumber}) and (\ref{self-intersection}) that 
\[
\begin{array}{llll}
\rho_{\Delta'} = 13-3 = 10, & 
D_1'^2 = 8, \, D_2'^2 = D_3'^2 = D_4'^2 = -2, \, D_5'^2 = 2, \\
 & D_6'^2 = \cdots = D_{13}'^2 =-2. 
\end{array}
\]
Let $L'$ be a lattice generated by divisors $\{D_i' \}_{i=1}^{13}$. 
By solving the equation (\ref{ToricDivisorsRelations}), one sees that $\{ D_4', D_5', D_3', D_{11}', D_{10}', D_9', D_6', D_{12}', D_{13}', D_2'  \}$ form a basis for $L'$. 
By taking a new basis
\[
\left\{
\begin{array}{l}
D_4',\, -D_3'+D_5',\, D_3'+D_4'+D_{11}',\\
-3D_3'-2D_4'+D_5'+D_6'-D_9'-2D_{10}'-3D_{11}'+D_{13}',\\
D_3'-D_5'+D_{10}',\\
3D_3'+2D_4'-D_5'-2 D_6'+D_{10}'+2D_{11}'-D_{12}'-D_{13}', \\
5D_3'+3D_4'-2D_5'-2D_6'+2D_9'+3D_{10}'+4D_{11}'-D_{12}'-D_{13}', \\
D_2'-5D_3'-3D_4'+2D_5'+2D_6'-D_9'-3D_{10}'-4D_{11}'+D_{12}'+D_{13}', \, D_{12}',\\
-D_2'+3D_3'+2D_4'-D_5'-D_6'+D_9'+2D_{10}'+3D_{11}'-D_{12}'-D_{13}'
\end{array}
\right\}, 
\]
one sees that the lattice $L'$ is isometric to $U\oplus A_1 \oplus E_7$, which is hyperbolic and a primitive sublattice of the $K3$ lattice. 
Therefore, $\Pic_{\!\Delta'}\simeq U\oplus A_1\oplus E_7$. 
By similar computation, one has $\Pic_{\!\Delta}\simeq U\oplus A_1\oplus E_7$. \\

\noindent
\underline{Case 2}\quad 
We have $\Delta\simeq\Delta'$. 

Set one-simplices of $\Sigma'$ as follows: 
\[
\begin{array}{lll}
m_1=( -1,  -1,  1), & m_2=( -1,  -1,  -1), & m_3=( 3,  -1,  -1), \\
m_4=( 3,  0,  -1), & m_5=( 0,  2,  -1), & m_6=( -1,  2,  -1), \\
m_7=( -1,  -1,  0), & m_8=( 0,  -1,  -1), & m_9=( 1,  -1,  -1), \\
m_{10}=( 2,  -1,  -1), & m_{11}=( 1,  -1,  0), & m_{12}=( -1,  0,  -1), \\
m_{13}=( -1,  1,  -1), 
\end{array}
\]
and let $\tilde{D_i'}$ be the toric divisor determined by the lattice point $m_i$ for $i=1,\ldots,13$, and $D_i':=\tilde{D_i'}|_{{-}K_X}$ with $X:=\widetilde{\mathbb{P}_{\Sigma'}}$. 
It can be easily seen by formulas (\ref{PicardNumber}) and (\ref{self-intersection}) that 
\[
\begin{array}{llll}
\rho_{\Delta'} = 13-3 = 10, & 
D_1'^2 = 8, \, D_2'^2 = D_3'^2 = -2, \, D_4'^2 = D_5'^2 = 0, \\
&  D_6'^2 =\cdots = D_{13}'^2 =-2. 
\end{array}
\]
Let $L'$ be a lattice generated by divisors $\{D_i' \}_{i=1}^{13}$. 
By solving the equation (\ref{ToricDivisorsRelations}), one sees that $\{ D_4', D_3', D_{11}', D_9', D_8', D_2', D_7', D_{12}', D_{13}', D_6'  \}$ form a basis for $L'$. 
By taking a new basis
\[
\left\{
\begin{array}{l}
D_4',\, D_3'+D_4',\, D_{11}'-D_4',\, D_9',\, D_8',\, D_2',\, D_7',\, D_{12}',D_{13}',\, D_6'
\end{array}
\right\}, 
\]
one sees that the lattice $L'$ is isometric to $U\oplus A_1 \oplus E_7$, which is hyperbolic and a primitive sublattice of the $K3$ lattice. 
Therefore, $\Pic_{\!\Delta'}\simeq U\oplus A_1\oplus E_7$. 
By similar computation, one has $\Pic_{\!\Delta'}\simeq U\oplus A_1\oplus E_7$. \\

\noindent
\underline{Case 3}\quad
Set one-simplices of $\Sigma$ as follows: 
\[
\begin{array}{lll}
v_1=( -1,  -1, 1), & v_2=( -1, -1,  -1),  & v_3=( 3,  -1,  -1), \\
v_4=( 3,  0,  -1), & v_5=( 0,  2, -1), & v_6=( -1,  1,  -1), \\
v_7=( -1,  -1,  0), & v_8=( 0,  -1,  -1), & v_9=( 1,  -1,  -1), \\
v_{10}=( 2,  -1,  -1), & v_{11}=( 1,  -1,  0), & v_{12}=( -1,  0,  -1), \\
v_{13}=( -1,  0,  0), 
\end{array}
\]
and let $\tilde{D_i}$ be the toric divisor determined by the lattice point $v_i$ for $i=1,\ldots,13$, and $D_i:=\tilde{D_i}|_{{-}K_X}$ with $X=\widetilde{\mathbb{P}_{\Sigma}}$. 
It can be easily seen by formulas (\ref{PicardNumber}) and (\ref{self-intersection}) that 
\[
\begin{array}{llll}
\rho_{\Delta} = 13-3 = 10, & 
D_1^2 = 8, \, D_2^2 = D_3^2 = -2, \, D_4^2 = 0, \, D_5^2 = 2, \\
&  D_6^2 = \cdots = D_{13}^2 =-2. 
\end{array}
\]
Let $L$ be a lattice generated by divisors $\{D_i \}_{i=1}^{13}$. 
By solving the equation (\ref{ToricDivisorsRelations}), one sees that $\{ D_4, D_3, D_{11}, D_9, D_8, D_2, D_7, D_{12}, D_6, D_{13}  \}$ form a basis for $L$. 
By taking a new basis
\[
\left\{
\begin{array}{l}
D_4,\, D_3+D_4,\, D_{11}-D_4,\, D_9,\, D_8,\, D_2,\, D_7,\, D_{12},\, D_6,\, D_{13}
\end{array}
\right\}, 
\]
one sees that the lattice $L$ is isometric to $U\oplus A_1 \oplus E_7$, which is hyperbolic and a primitive sublattice of the $K3$ lattice. 
Therefore, $\Pic_{\!\Delta}\simeq U\oplus A_1\oplus E_7$. \\

Set one-simplices of $\Sigma'$ as follows: 
\[
\begin{array}{lll}
m_1=( -1,  -1,  1), & m_2=( -1,  -1,  -1), & m_3=( 4,  -1,  -1), \\
m_4=( 3,  0,  -1), & m_5=( 0,  2,  -1), & m_6=( -1,  2,  -1), \\
m_7=( -1,  -1,  0), & m_8=( 0,  -1,  -1), & m_9=( 1,  -1,  -1), \\
m_{10}=( 2,  -1,  -1), & m_{11}=( 3,  -1,  -1), & m_{12}=( -1,  0,  -1), \\
m_{13}=( -1,  1,  -1), 
\end{array}
\]
and let $\tilde{D_i'}$ be the toric divisor determined by the lattice point $m_i$ for $i=1,\ldots,13$, and $D_i':=\tilde{D_i'}|_{{-}K_X}$ with $X:=\widetilde{\mathbb{P}_{\Sigma'}}$. 
It can be easily seen by formulas (\ref{PicardNumber}) and (\ref{self-intersection}) that 
\[
\begin{array}{llll}
\rho_{\Delta'} = 13-3 = 10, & 
D_1'^2 = 8, \, D_2'^2 = D_3'^2 = D_4'^2 = -2, \, D_5'^2 = 0, \\
&  D_6'^2 = \cdots = D_{13}'^2 =-2. 
\end{array}
\]
Let $L'$ be a lattice generated by divisors $\{D_i' \}_{i=1}^{13}$. 
By solving the equation (\ref{ToricDivisorsRelations}), one sees that $\{ D_4', D_3', D_{11}', D_9', D_8', D_2', D_7', D_{12}', D_{13}', D_6'  \}$ form a basis for $L'$. 
By taking a new basis
\[
\left\{
\begin{array}{l}
D_3'+D_4'+D_{11}',\, D_3'+D_4',\, D_4',\, D_9',\, D_8',\, D_2',\, D_7',\, D_{12}',\, D_{13}',\, D_6'
\end{array}
\right\}, 
\]
one sees that the lattice $L'$ is isometric to $U\oplus A_1 \oplus E_7$, which is hyperbolic and a primitive sublattice of the $K3$ lattice. 
Thus, $\Pic_{\!\Delta'}\simeq U\oplus A_1\oplus E_7$. 

It is well-known that the lattice $U\oplus A_1\oplus E_7$ is a primitive sublattice of the $K3$ lattice $\Lambda_{K3}$. 
Moreover, by Lemma~\ref{Nishiyama}, the relation 
$(\Pic_{\!\Delta})^\perp_{\Lambda_{K3}}\simeq U^\perp_{U^{\oplus3}}\oplus(A_1)^\perp_{E_8}\oplus (E_7)^\perp_{E_8}=U^{\oplus2}\oplus E_7\oplus A_1\simeq U\oplus\Pic_{\!\Delta'}$ holds. 

\subsection{No. 26}
In all cases, we set one-simplices of fans in terms of a basis of $M_{(1,3,4,5)}\otimes\mathbb{R}$
\[
(-3,1,0,0),\, (-4,0,1,0),\, (-5,0,0,1). 
\]

\begin{lem}\label{Lemma26}
If lattices $L$ and $L'$ have the signature, the discriminant, and the rank of $L$ and $L'$ are respectively $(1,\, 9)$, $\discr{L}=\discr{L'}=-13$, and $\rk{L}=\rk{L'}=10$, then, the lattices are primitive sublattices of the $K3$ lattice and $U\oplus L'$ is the orthogonal complement of $L$. 
\end{lem}
\proof
Note that the discriminant groups $A_L,\, A_{L'}$ of $L$ and $L'$ are isomorphic to $\mathbb{Z}\slash13\mathbb{Z}$, and that the minimal number of the generators is $l(A_L) = l(A_{L'}) =1$. 
Since the signature of $L$ and $L'$ is $(t_+,\, t_-)=(1,9)$ and the rank is $\rk{L}=\rk{L'}=10$, we have 
\begin{eqnarray*}
19-t_-=10\geq 0, \quad 3-t_+=2\geq 0, \quad \textnormal{and}\\
22-\rk{L}=22-\rk{L'}=12>1=l(A_L)=l(A_{L'}),
\end{eqnarray*}
by Corollary~\ref{NikulinPrimitive}, the first statement is shown. 
Since the discriminant of $U$ is $-1$, we have $\discr{(U\oplus L')} = -\discr{L'}=13=-\discr{L}$, and thus by Corollary~\ref{NikulinOrthogonal}, the last assertion is proved. 
\QED

\noindent
\underline{Case 1}\quad 
Set one-simplices of $\Sigma$ as follows: 
\[
\begin{array}{lll}
v_1=( -1,  -1,  1), & v_2=( -1,  -1,  -1), & v_3=( 2,  -1,  -1), \\
v_4=( 2,  0,  -1), & v_5=( -1,  2,  -1), & v_6=( -1,  1, 0), \\
v_7=( 0,  -1,  1), & v_8=( -1,  -1,  0), & v_9=( 0,  -1,  -1), \\
v_{10}=( 1,  -1,  -1), & v_{11}=( -1,  0,  -1), & v_{12}=( -1,  1,  -1), \\
v_{13}=( 1,  -1,  0), 
\end{array}
\]
and let $\tilde{D_i}$ be the toric divisor determined by the lattice point $v_i$ for $i=1,\ldots,13$, and $D_i:=\tilde{D_i}|_{{-}K_X}$ with $X:=\widetilde{\mathbb{P}_{\Sigma}}$. 
It can be easily seen by formulas (\ref{PicardNumber}) and (\ref{self-intersection}) that 
\[
\begin{array}{llll}
\rho_{\Delta} = 13-3 = 10, &
D_1^2 = D_2^2 = D_3^2 = -2,  \, D_4^2 = 0, \, D_5^2 = D_6^2 = -2, \\
&  D_7^2 = 0, \, D_8^2 = \cdots = D_{13}^2 =-2. 
\end{array}
\]
Let $L$ be a lattice generated by divisors $\{D_i \}_{i=1}^{13}$. 
By solving the equation (\ref{ToricDivisorsRelations}), one sees that $\{$ $D_7$, $D_4$, $D_3$, $D_{13}$, $D_{10}$, $D_9$, $D_2$, $D_8$, $D_{11}$, $D_{12}$ $\}$ form a basis for $L$. 
By taking a new basis
\[
\left\{
\begin{array}{l}
D_7,\, D_7+D_{13},\, D_3-D_7,\, D_3-D_4+3D_7+D_9+D_{10}+2D_{13},\\
D_{10},\, D_2+D_8+D_9,\, D_2,\, D_8,\, D_{11},\, D_{12}
\end{array}
\right\}, 
\]
one sees that the lattice $L$ is isometric to $U\oplus \tilde{L}$ with some lattice $\tilde{L}$. 
By a direct computation, one sees that $\sign{L}=(1,\, 9),\, \discr{L}=-13$, and $\rk{L}=10$, thus, $\discr{\tilde{L}}=13$ and $\rk{\tilde{L}}=8$ hold. 
In particular, the discriminant group $A_L$ of $L$ is isomorphic to $\mathbb{Z}\slash13\mathbb{Z}$, and $l(A_L)=1$. 

Set one-simplices of $\Sigma'$ as follows: 
\[
\begin{array}{lll}
m_1=(-1,  -1,  1), & m_2=( -1,  -1,  -1), & m_3=( 2,  -1,  -1), \\
m_4=( 2,  0,  -1), & m_5=( -1,  1,  0), & m_6=( 0,  -1,  1), \\
m_7=( -1,  0,  -1), & m_8=( -1,  -1,  0), & m_9=( 0,  -1,  -1), \\
m_{10}=( 1,  -1,  -1), & m_{11}=( 1,  0,  -1), & m_{12}=( 0,  0,  -1), \\
m_{13}=( 1,  -1,  0), 
\end{array}
\]
and let $\tilde{D_i'}$ be the toric divisor determined by the lattice point $m_i$ for $i=1,\ldots,13$, and $D_i':=\tilde{D_i'}|_{{-}K_X}$ with $X:=\widetilde{\mathbb{P}_{\Sigma'}}$. 
It can be easily seen by formulas (\ref{PicardNumber}) and (\ref{self-intersection}) that 
\[
\begin{array}{llll}
\rho_{\Delta'} = 13-3 = 10, & 
D_1'^2 = D_2'^2 = D_3'^2 = -2, \, D_4'^2 = 0, \, 
D_5'^2 = 4, \, D_6'^2 = 0, \\
& D_7'^2 = \cdots = D_{13}'^2 =-2. 
\end{array}
\]
Let $L'$ be a lattice generated by divisors $\{D_i' \}_{i=1}^{13}$. 
By solving the equation (\ref{ToricDivisorsRelations}), one sees that $\{$ $D_4'$, $D_6'$, $D_{13}'$, $D_3'$, $D_{10}'$, $D_9'$, $D_2'$, $D_8'$, $D_7'$, $D_{12}'$ $\}$ form a basis for $L'$. 
By taking a new basis
\[
\left\{
\begin{array}{l}
D_4',\, D_3'+D_4',\, D_4'-D_{13}',\, 2D_3'+2D_4'-D_6'+D_{10}'+D_{13}',\\
-D_4'+D_{10}',\, D_9',\, D_2',\, D_8',\, D_7',\, D_{12}'
\end{array}
\right\}, 
\]
one sees that the lattice $L'$ is isometric to $U\oplus \tilde{L'}$ with some lattice $\tilde{L'}$. 
By a direct computation, one sees that $\sign{L'}=(1,\, 9),\, \discr{L'}=-13$, and $\rk{L'}=10$, thus, $\discr{\tilde{L'}}=13$ and $\rk{\tilde{L'}}=8$ hold. 
In particular, the discriminant group $A_{L'}$ of $L'$ is isomorphic to $\mathbb{Z}\slash13\mathbb{Z}$, and $l(A_{L'})=1$. 

\noindent
\underline{Case 2}\quad
Set one-simplices of $\Sigma$ as follows: 
\[
\begin{array}{lll}
v_1=( -1,  -1,  0), & v_2=( -1,  -1,  -1), & v_3=( 3,  -1,  -1), \\
v_4=( 2,  0,  -1), & v_5=( -1,  2,  -1), & v_6=( -1,  1,  0), \\
v_7=( 0,  -1,  1), & v_8=( -1,  0,  0), & v_9=( 0,  -1,  -1), \\
v_{10}=( 1,  -1,  -1), & v_{11}=( 2,  -1,  -1), & v_{12}=( -1,  0,  -1), \\
v_{13}=( -1,  1,  -1), 
\end{array}
\]
and let $\tilde{D_i}$ be the toric divisor determined by the lattice point $v_i$ for $i=1,\ldots,13$, and $D_i:=\tilde{D_i}|_{{-}K_X}$ with $X:=\widetilde{\mathbb{P}_{\Sigma}}$. 
It can be easily seen by formulas (\ref{PicardNumber}) and (\ref{self-intersection}) that 
\[
\rho_{\Delta} = 13-3 = 10, \quad
D_1^2 =\cdots = D_6^2 = -2, \, D_7^2 = 2, \, D_8^2 = \cdots = D_{13}^2 =-2. 
\]
Let $L$ be a lattice generated by divisors $\{D_i \}_{i=1}^{13}$. 
By solving the equation (\ref{ToricDivisorsRelations}), one sees that $\{$ $D_7$, $D_1$, $D_8$, $D_2$, $D_{12}$, $D_{13}$, $D_9$, $D_{10}$, $D_{11}$, $D_3$ $\}$ form a basis for $L$. 
By taking a new basis
\[
\left\{
\begin{array}{l}
D_7-D_8,\, D_3+D_7-D_8,\, D_1,\, D_2,\, D_{12},\, D_{13},\, D_9,\\
D_{10},\, -D_7+D_8+D_{11},\, -2D_3-4D_7+5D_8
\end{array}
\right\}, 
\]
one sees that the lattice $L$ is isometric to $U\oplus \tilde{L}$ with some lattice $\tilde{L}$. 
By a direct computation, one sees that $\sign{L}=(1,\, 9),\, \discr{L}=-13$, and $\rk{L}=10$, thus, $\discr{\tilde{L}}=13$ and $\rk{\tilde{L}}=8$ hold. 
In particular, the discriminant group $A_L$ of $L$ is isomorphic to $\mathbb{Z}\slash13\mathbb{Z}$, and $l(A_L)=1$. 

Set one-simplices of $\Sigma'$ as follows: 
\[
\begin{array}{lll}
m_1=(-1,  -1,  0), & m_2=( -1,  -1, -1), & m_3=( 3,  -1,  -1), \\
m_4=( 2, 0,  -1), & m_5=( -1,  1,  0), & m_6=( 0,  -1,  1), \\
m_7=( -1,  0,  -1), & m_8=( -1,  0,  0), & m_9=( 0,  -1,  -1), \\
m_{10}=( 1,  -1,  -1), & m_{11}=( 2,  -1,  -1), & m_{12}=( 0,  0,  -1), \\
m_{13}=( 1,  0,  -1), 
\end{array}
\]
and let $\tilde{D_i'}$ be the toric divisor determined by the lattice point $m_i$ for $i=1,\ldots,13$, and $D_i':=\tilde{D_i'}|_{{-}K_X}$ with $X:=\widetilde{\mathbb{P}_{\Sigma'}}$. 
It can be easily seen by formulas (\ref{PicardNumber}) and (\ref{self-intersection}) that 
\[
\begin{array}{llll}
\rho_{\Delta'} = 13-3 = 10, & 
D_1'^2 = \cdots = D_4'^2 = -2, \, D_5'^2 = 4, \, D_6'^2 = 2, \\
 &  D_7'^2 = \cdots = D_{13}'^2 =-2. 
\end{array}
\]
Let $L'$ be a lattice generated by divisors $\{D_i' \}_{i=1}^{13}$. 
By solving the equation (\ref{ToricDivisorsRelations}), one sees that $\{$ $D_6'$, $D_1'$, $D_8'$, $D_4'$, $D_3'$, $D_{11}'$, $D_{10}'$, $D_9'$, $D_{13}'$, $D_{12}'$ $\}$ form a basis for $L'$. 
By taking a new basis
\[
\left\{
\begin{array}{l}
D_3'+D_4',\, D_3'+D_4'+D_{11}',\, -4D_3'-4D_4'+D_6'-2D_{11}',\\
-D_4',\,  D_1',\,  D_8',\,  -D_3'-D_4'+D_{10}',\,  D_9',\\
 -2D_3'-2D_4'-D_{11}'+D_{13}',\,  D_{12}' 
\end{array}
\right\}, 
\]
one sees that the lattice $L'$ is isometric to $U\oplus \tilde{L'}$ with some lattice $\tilde{L'}$. 
By a direct computation, one sees that $\sign{L'}=(1,\, 9),\, \discr{L'}=-13$, and $\rk{L'}=10$, thus, $\discr{\tilde{L'}}=13$ and $\rk{\tilde{L'}}=8$ hold. 
In particular, the discriminant group $A_{L'}$ of $L'$ is isomorphic to $\mathbb{Z}\slash13\mathbb{Z}$, and $l(A_{L'})=1$. 

\noindent
\underline{Case 3}\quad 
Set one-simplices of $\Sigma$ as follows: 
\[
\begin{array}{lll}
v_1=( -1,  -1,  1), & v_2=( -1,  -1,  -1), & v_3=( 3,  -1,  -1), \\
v_4=( 2,  0,  -1), & v_5=( -1,  1,  0), & v_6=( 0,  -1,  1), \\
v_7=( -1,  0,  -1), & v_8=( -1,  -1,  0), & v_9=( 0,  -1,  -1), \\
v_{10}=( 1,  -1,  -1), & v_{11}=( 2,  -1,  -1), & v_{12}=( 0,  0,  -1), \\
v_{13}=( 1, 0,  -1), 
\end{array}
\]
and let $\tilde{D_i}$ be the toric divisor determined by the lattice point $v_i$ for $i=1,\ldots,13$, and $D_i:=\tilde{D_i}|_{{-}K_X}$ with $X:=\widetilde{\mathbb{P}_{\Sigma}}$. 
It can be easily seen by formulas (\ref{PicardNumber}) and (\ref{self-intersection}) that 
\[
\begin{array}{llll}
\rho_{\Delta} = 13-3 = 10, & 
D_1^2 = \cdots = D_4^2 = -2, \, D_5^2 = 4, \, D_6^2 = 0, \\
 & D_7^2 = \cdots = D_{13}^2 =-2. 
\end{array}
\]
Let $L$ be a lattice generated by divisors $\{D_i \}_{i=1}^{13}$. 
By solving the equation (\ref{ToricDivisorsRelations}), one sees that $\{$  $D_8$, $D_1$, $D_6$, $D_4$, $D_{13}$, $D_{12}$, $D_3$, $D_{11}$, $D_{10}$, $D_9$ $\}$ form a basis for $L$. 
By taking a new basis
\[
\left\{
\begin{array}{l}
D_6,\, D_4+D_6,\, D_1-D_3-D_4,\, D_8,\, -D_6+D_{13},\\
D_{12},\, D_3-D_4-4D_6,\, D_{11},\, D_{10},\, D_9
\end{array}
\right\}, 
\]
one sees that the lattice $L$ is isometric to $U\oplus \tilde{L}$ with some lattice $\tilde{L}$. 
By a direct computation, one sees that $\sign{L}=(1,\, 9),\, \discr{L}=-13$, and $\rk{L}=10$, thus, $\discr{\tilde{L}}=13$ and $\rk{\tilde{L}}=8$ hold. 
In particular, the discriminant group $A_L$ of $L$ is isomorphic to $\mathbb{Z}\slash13\mathbb{Z}$, and $l(A_L)=1$. 

Set one-simplices of $\Sigma'$ as follows: 
\[
\begin{array}{lll}
m_1=(-1,  -1, 0), & m_2=( -1,  -1,  -1), & m_3=( 2,  -1,  -1), \\
m_4=( 2, 0,  -1), & m_5=( -1,  2,  -1), & m_6=( -1,  1,  0), \\
m_7=( 0,  -1,  1), & m_8=( -1,  0,  0), & m_9=( 0,  -1,  -1), \\
m_{10}=( 1,  -1,  -1), & m_{11}=( 1,  -1,  0), & m_{12}=( -1,  0,  -1), \\
m_{13} =( -1,  1,  -1), 
\end{array}
\]
and let $\tilde{D_i'}$ be the toric divisor determined by the lattice point $m_i$ for $i=1,\ldots,13$, and $D_i':=\tilde{D_i'}|_{{-}K_X}$ with $X:=\widetilde{\mathbb{P}_{\Sigma'}}$. 
It can be easily seen by formulas (\ref{PicardNumber}) and (\ref{self-intersection}) that 
\[
\begin{array}{llll}
\rho_{\Delta'} = 13-3 = 10, & 
D_1'^2 = D_2'^2 = D_3'^2 = -2, \, D_4'^2 = 0, \, D_5'^2 = D_6'^2 = -2, \\
&  D_7'^2 = 2, \, D_8'^2 =\cdots = D_{13}'^2 =-2. 
\end{array}
\]
Let $L'$ be a lattice generated by divisors $\{D_i' \}_{i=1}^{13}$. 
By solving the equation (\ref{ToricDivisorsRelations}), one sees that $\{$  $D_7'$, $D_4'$, $D_3'$, $D_{11}'$, $D_{10}'$, $D_9'$, $D_2'$, $D_1'$, $D_{12}'$, $D_{13}'$ $\}$ form a basis for $L'$. 
By taking a new basis
\[
\left\{
\begin{array}{l}
D_3'+D_4',\, D_4',\, -2D_3'-4D_4'+D_7',\, -D_4'+D_{11}',\\
-D_4'+D_{10}',\, D_9',\, D_2',\, D_1',\, D_{12}',\, D_{13}'
\end{array}
\right\}, 
\]
one sees that the lattice $L'$ is isometric to $U\oplus \tilde{L'}$ with some lattice $\tilde{L'}$. 
By a direct computation, one sees that $\sign{L'}=(1,\, 9),\, \discr{L'}=-13$, and $\rk{L'}=10$, thus, $\discr{\tilde{L'}}=13$ and $\rk{\tilde{L'}}=8$ hold. 
In particular, the discriminant group $A_{L'}$ of $L'$ is isomorphic to $\mathbb{Z}\slash13\mathbb{Z}$, and $l(A_{L'})=1$. 

\noindent
\underline{Case 4}\quad
Set one-simplices of $\Sigma$ as follows: 
\[
\begin{array}{lll}
v_1=( -1,  -1,  1), & v_2=( -1,  -1,  -1), & v_3=( 3,  -1,  -1), \\
v_4=( 2,  0,  -1), & v_5=( -1,  2,  -1), & v_6=( -1,  1,  0), \\
v_7=( 0,  -1,  1), & v_8=( -1,  -1,  0), & v_9=( 0,  -1,  -1), \\
v_{10}=( 1,  -1,  -1), & v_{11}=( 2,  -1,  -1), & v_{12}=( -1,  0,  -1), \\
v_{13}=( -1,  1,  -1), 
\end{array}
\]
and let $\tilde{D_i}$ be the toric divisor determined by the lattice point $v_i$ for $i=1,\ldots,13$, and $D_i:=\tilde{D_i}|_{{-}K_X}$ with $X:=\widetilde{\mathbb{P}_{\Sigma}}$. 
One can easily seen by formulas (\ref{PicardNumber}) and (\ref{self-intersection}) that 
\[
\begin{array}{llll}
\rho_{\Delta} = 13-3 = 10, & 
D_1 = \cdots = D_6 = -2, \, D_7 = 0, \\
&  D_8 =\cdots = D_{13} =-2. 
\end{array}
\]
Let $L$ be a lattice generated by divisors $\{D_i \}_{i=1}^{13}$. 
By solving the equation (\ref{ToricDivisorsRelations}), one sees that $\{$ $D_7$, $D_4$, $D_3$, $D_{11}$, $D_{10}$, $D_9$, $D_2$, $D_{12}$, $D_{13}$, $D_8$  $\}$ form a basis for $L$. 
By taking a new basis
\[
\left\{
\begin{array}{l}
D_7, \,D_4+D_7,\, D_3-D_4-4D_7,\, D_{11},\, D_{10},\, D_9,\, D_2,\, D_{12},\, D_{13},\, D_8
\end{array}
\right\}, 
\]
one sees that the lattice $L$ is isometric to $U\oplus\tilde{L}$ with some lattice $\tilde{L}$. 
By a direct computation, one sees that $\sign{L}=(1,\, 9),\, \discr{L}=-13$, and $\rk{L}=10$, and thus, $\discr{\tilde{L}}=13$ and $\rk{\tilde{L}}=8$ hold. 
In particular, the discriminant group $A_L$ of $L$ is isomorphic to $\mathbb{Z}\slash13\mathbb{Z}$, and $l(A_L)=1$. 

Set one-simplices of $\Sigma'$ as follows: 
\[
\begin{array}{lll}
m_1=( -1,  -1,  0), & m_2=( -1,  -1,  -1), & m_3=( 2,  -1,  -1), \\
m_4=( 2,  0,  -1), & m_5=( -1,  1,  0), & m_6=( 0,  -1,  1), \\
m_7=( -1,  0,  -1), & m_8=( 0,  -1,  -1), & m_9 =( 1,  -1,  -1), \\
m_{10}=( -1,  0,  0), & m_{11}=( 1,  -1,  0), & m_{12}=( 0,  0,  -1), \\
m_{13}=( 1,  0,  -1), 
\end{array}
\]
and let $\tilde{D_i'}$ be the toric divisor determined by the lattice point $m_i$ for $i=1,\ldots,13$, and $D_i':=\tilde{D_i'}|_{{-}K_X}$ with $X:=\widetilde{\mathbb{P}_{\Sigma'}}$. 
It can be easily seen by formulas (\ref{PicardNumber}) and (\ref{self-intersection}) that 
\[
\begin{array}{llll}
\rho_{\Delta'} = 13-3 = 10, & 
D_1'^2 = D_2'^2 = D_3'^2 = -2, \, D_4'^2 = 0, \, D_5'^2 = 4, \\
&  D_6'^2 = 2, \, D_7'^2 = \cdots = D_{13}'^2 =-2. 
\end{array}
\]
Let $L'$ be a lattice generated by divisors $\{D_i' \}_{i=1}^{13}$. 
By solving the equation (\ref{ToricDivisorsRelations}), one sees that $\{$  $D_{11}'$, $D_3'$, $D_4'$, $D_{13}'$, $D_{12}'$, $D_7'$, $D_2'$, $D_8'$, $D_1'$, $D_{10}'$  $\}$ form a basis for $L'$. 
By taking a new basis
\[
\left\{
\begin{array}{l}
D_4',\, D_4'+D_{13}',\, D_{11}',\, D_3'-2D_4'-D_{13}',\\
-D_4'+D_{12}',\, D_7',\, D_2',\, D_8',\, D_1',\, D_{10}'
\end{array}
\right\}, 
\]
one sees that the lattice $L'$ is isometric to $U\oplus \tilde{L'}$ with some lattice $\tilde{L'}$. 
By a direct computation, one sees that $\sign{L'}=(1,\, 9),\, \discr{L'}=-13$, and $\rk{L'}=10$, thus, $\discr{\tilde{L'}}=13$ and $\rk{\tilde{L'}}=8$ hold. 
In particular, the discriminant group $A_{L'}$ of $L'$ is isomorphic to $\mathbb{Z}\slash13\mathbb{Z}$, and $l(A_{L'})=1$. 

In all cases $1$ to $4$, we obtain lattices $L$ and $L'$ satisfying assumptions in Lemma~\ref{Lemma26}. 
Therefore, we can conclude that $\Pic_{\!\Delta}=U\oplus\tilde{L},\,  \Pic_{\!\Delta'}=U\oplus\tilde{L'}$, with $\discr{\tilde{L}}=\discr{\tilde{L'}}=13$ and $\rk{\tilde{L}}=\rk{\tilde{L'}}=8$, and that the relation $(\Pic_{\!\Delta})^\perp_{\Lambda_{K3}}\simeq U\oplus \Pic_{\!\Delta'}$ holds. 

\subsection{Nos. 35, 36, and 37}
Set one-simplices of $\Sigma$ in terms of a basis of $M_{(1,1,4,6)}\otimes\mathbb{R}$ 
\[
(-1,1,0,0),\, (-4,0,1,0),\, (-6,0,0,1) :
\]
\[
\begin{array}{lll}
v_1 = ( -1,  -1,  1), & v_2 = ( -1,  -1,  -1), & v_3 = ( 11,  -1,  -1), \\
v_4 = ( -1,  2,  -1), & v_5 = ( -1,  -1,  0), & v_6 = ( 5,  -1,  0), \\
v_7 = ( 0,  -1,  -1), & v_8 = ( 1,  -1, -1), & v_9 = ( 2,  -1,  -1), \\
v_{10} = ( 3,  -1,  -1), & v_{11} = ( 4,  -1,  -1), & v_{12} = ( 5,  -1,  -1), \\
v_{13} = ( 6,  -1, -1), & v_{14} = ( 7,  -1,  -1),& v_{15} = ( 8,  -1,  -1), \\
v_{16} = ( 9,  -1,  -1), & v_{17} = ( 10,  -1,  -1), & v_{18} = ( 7,  0,  -1), \\
v_{19} = ( 3,  1,  -1), & v_{20} = ( -1,  0,  -1), & v_{21} = ( -1,  1,  -1), 
\end{array}
\]
and let $\tilde{D_i}$ be the toric divisor determined by the lattice point $v_i$ for $i=1,\ldots,21$, and $D_i:=\tilde{D_i}|_{{-}K_X}$ with $X:=\widetilde{\mathbb{P}_{\Sigma}}$. 
It can be easily seen by formulas (\ref{PicardNumber}) and (\ref{self-intersection}) that 
\[
\begin{array}{llll}
\rho_{\Delta} = 21-3 = 18, & 
D_1^2 = 2, \, D_2^2 = D_3^2 = -2 , \, D_4^2 = 0 \\
&  D_5^2 = \cdots = D_{21}^2 = -2. 
\end{array}
\]
Let $L$ be a lattice generated by divisors $\{D_i \}_{i=1}^{21}$. 
By solving the equation (\ref{ToricDivisorsRelations}), one sees that $\{$ $D_4$, $D_1$, $D_5$,$D_{21}$,$D_{20}$,$D_{19}$,$D_{18}$,$D_3$,$D_{17}$,$D_{16}$,$D_{15}$,$D_{14}$,$D_{13}$, \\ $D_{12}$,$D_{11}$,$D_{10}$,$D_9$,$D_8$ $\}$ form a basis for $L$. 
By taking a new basis
\[
\left\{
\begin{array}{l}
D_4,\, D_4+D_{21},\, -D_5,\, -D_1+2D_4+D_{19}+D_{21},\, -D_4+D_{20},\\
D_3-D_4+D_{16}+D_{17}+D_{18}+D_{19}-D_{20}-D_{21},\, D_{18},\, D_3,\\ 
D_{17},\, D_{14}+D_{15}+D_{16},\, D_{15},\, D_{14},\, D_{13}+D_{14}+D_{15},\\
D_{12},\, D_{11},\, D_{10},\, D_9,\, D_8
\end{array}
\right\}, 
\]
one sees that the lattice $L$ is isometric to $U\oplus \tilde{L}$, where $\tilde{L}$ is a negative-definite of rank $16$ and discriminant $1$. 
By the classification of unimodular lattices, we have $\tilde{L}\simeq E_8^{\oplus 2}$. 
Therefore, $\Pic_{\!\Delta}\simeq U\oplus E_8^{\oplus 2}$. 

Set one-simplices of $\Sigma'$ in terms of a basis of $M_{(3,5,11,14)}\otimes\mathbb{R}$ 
\[
(1,0,1,-1),  \, (2,1,-1,0), \, (10,-1,-1,-1) : 
\]
\[
\begin{array}{lll}
m_1 = (-1, 0, 0), & m_2 = (0, 0, 1), & m_3 = (2, 4, -1), \\
m_4 = (1, -1, 0), & m_5 = (1, 2, 0), 
\end{array}
\]
and let $\tilde{D_i'}$ be the toric divisor determined by the lattice point $m_i$ for $i=1,\ldots,5$, and $D_i':=\tilde{D_i'}|_{{-}K_X}$ with $X:=\widetilde{\mathbb{P}_{\Sigma'}}$. 
It can be easily seen by formulas (\ref{PicardNumber}) and (\ref{self-intersection}) that 
\[
\rho_{\Delta'} = 5-3 = 2, \qquad
D_1'^2 = 18, \, D_2'^2 = D_3'^2 = 0, \, D_4'^2 = 8, \, D_5'^2 = -2. 
\]
Let $L'$ be a lattice generated by divisors $\{D_i' \}_{i=1}^{5}$. 
By solving the equation (\ref{ToricDivisorsRelations}), one sees that $\{$ $D_3'$, $D_5'$ $\}$ form a basis for $L'$. 
By taking a new basis $\left\{D_3', \, D_3'+D_5'\right\}$, one sees that the lattice $L'$ is isometric to $U$, which is a hyperbolic primitive sublattice of the $K3$ lattice. 
Thus, $\Pic_{\!\Delta'}\simeq U$. 

It is well-known that lattices $U$ and $U\oplus E_8^{\oplus 2}$ are primitive sublattices of the $K3$ lattice $\Lambda_{K3}$ and it is clear that the relation 
$(\Pic_{\!\Delta})^\perp_{\Lambda_{K3}}\simeq U^\perp_{U^{\oplus3}}\oplus(E_8^{\oplus 2})^\perp_{E_8^{\oplus 3}}=U^{\oplus2}\simeq U\oplus\Pic_{\!\Delta'}$ holds. 

\subsection{Nos. 38 and 40}
Take bases of $M_{(1,1,3,5)}\otimes\mathbb{R}$, and of $M_{(3,4,10,13)}\otimes\mathbb{R}$, respectively: 
\begin{eqnarray*}
\{ (-1,1,0,0), \, (-3,0,1,0), \, (-5,0,0,1) \}, \\
\{ (1,0,1,-1), \, (3,1,0,-1), \, (9,-1,-1,-1)\}. 
\end{eqnarray*}

\begin{lem}\label{Lemma38-40}
If $L$ is a negative-definite lattice of rank $15$ of discriminant $-2$, then, it is a primitive sublattice of the $K3$ lattice. 
\end{lem}
\proof
Note that the discriminant group of $L$ is isomorphic to $\mathbb{Z}\slash2\mathbb{Z}$ of number of generator $l(A_L)=1$. 
Since the signature of $L$ and $L'$ is $(t_+,\, t_-)=(0,15)$ and the rank is $\rk{L}=15$, we have 
\[
19-t_-=4\geq 0, \quad 3-t_+=3\geq 0, \quad \textnormal{and}\quad
22-\rk{L}=7>1=l(A_L),
\]
by Corollary~\ref{NikulinPrimitive}, the assertion is proved. 
\QED

\noindent
\underline{Case 1} \quad
Set one-simplices of $\Sigma$ as follows: 
\[
\begin{array}{lll}
v_1=( -1,  -1,  1), & v_2=( -1,  -1,  -1), & v_3=( 9,  -1,  -1), \\
v_4=( 0,  2,  -1), & v_5=( -1,  2,  -1), & v_6=( -1,  -1,  0), \\
v_7=( 0,  -1,  -1), & v_8=( 1,  -1,  -1), & v_9=( 2,  -1,  -1), \\
v_{10}=( 3,  -1,  -1), & v_{11}=( 4,  -1,  -1), & v_{12}=( 5,  -1,  -1), \\
v_{13}=( 6,  -1,  -1), & v_{14}=( 7,  -1,  -1), & v_{15}=( 8,  -1,  -1), \\
v_{16}=( 6,  0,  -1), & v_{17}=( 3,  1,  -1), & v_{18}=( -1,  0,  -1), \\
v_{19}=( -1,  1, -1), & v_{20}=( 4,  -1,  0), 
\end{array}
\]
and let $\tilde{D_i}$ be the toric divisor determined by the lattice point $v_i$ for $i=1,\ldots,20$, and $D_i:=\tilde{D_i}|_{{-}K_X}$ with $X:=\widetilde{\mathbb{P}_{\Sigma}}$. 
It can be easily seen by formulas (\ref{PicardNumber}) and (\ref{self-intersection}) that 
\[
\rho_{\Delta} = 20-3 = 17, \qquad
D_1^2 = 2, \quad D_2^2 = \cdots = D_{20}^2 =-2. 
\]
Let $L$ be a lattice generated by divisors $\{D_i \}_{i=1}^{20}$. 
By solving the equation (\ref{ToricDivisorsRelations}), one sees that $\{$ $D_4$, $D_5$, $D_{19}$, $D_{18}$, $D_2$, $D_1$, $D_{20}$, $D_3$, $D_{16}$, $D_{15}$, $D_{14}$, $D_{13}$, $D_{12}$, $D_{11}$, $D_{10}$, $D_9$, $D_8$ $\}$ form a basis for $L$. 
By taking a new basis
\[
\left\{
\begin{array}{l}
D_4+D_5+D_{19},\, D_4+D_5,\, D_4,\, D_1-4D_4-4D_5-2D_{19},\\ 
D_2,\, -D_4-D_5+D_{18},\, D_{20},\, D_3,\, D_{16},\, D_{15},\, D_{14},\, D_{13},\\
D_{12},\, D_{11},\, D_{10},\, D_9,\, D_8
\end{array}
\right\}, 
\]
one sees that the lattice $L$ is isometric to $U\oplus \tilde{L}$, where $\discr{\tilde{L}}=-2$ and $\rk{\tilde{L}}=15$. 

Set one-simplices of $\Sigma'$ as follows: 
\[
\begin{array}{lll}
m_1=(1,0,0), & m_2=(0, 1, 0), & m_3=(0, 0,1), \\
m_4=(0, -2,  -3), & m_5=(-1,-3,-5), & m_6=(0,-1,-1),
\end{array}
\]
and let $\tilde{D_i'}$ be the toric divisor determined by the lattice point $m_i$ for $i=1,\ldots,6$, and $D_i':=\tilde{D_i'}|_{{-}K_X}$ with $X:=\widetilde{\mathbb{P}_{\Sigma'}}$. 
It can be easily seen by formulas (\ref{PicardNumber}) and (\ref{self-intersection}) that 
\[
\begin{array}{llll}
\rho_{\Delta} = 6-3 = 3, & 
D_1'^2 = 0, \, D_2'^2 = 6, \,  D_3'^2 = 16, \\
& D_4'^2 = -2, \, D_5'^2 = 0, \,  D_6'^2 =-2. 
\end{array}
\]
Let $L'$ be a lattice generated by divisors $\{D_i' \}_{i=1}^{6}$. 
By solving the equation (\ref{ToricDivisorsRelations}), one sees that $\{$ $D_1'$, $D_4'$, $D_6'$ $\}$ form a basis for $L'$. 
By taking a new basis $\left\{D_1',\, D_1'+D_4',\, D_6'-D_1' \right\}$, one sees that the lattice $L'$ is isometric to $U\oplus A_1$, which is hyperbolic and a primitive sublattice of the $K3$ lattice. 
Therefore, $\Pic_{\!\Delta'}\simeq U\oplus A_1$. 
Note that the discriminant group of $\Pic_{\!\Delta'}$ is isomorphic to  $\mathbb{Z}\slash 2\mathbb{Z}$ since $2$ is a prime number. \\

\noindent
\underline{Case 2}\quad
Set one-simplices of $\Sigma$ as follows: 
\[
\begin{array}{lll}
v_1=( -1,  -1,  -1),  & v_2=( 1,  1,  3),  & v_3=( 1,  3,  9),  \\
v_4=( 1,  3, -1),  & v_5=( 1,  0,  -1),  & v_6=( 0,  0,  1),  \\
v_7=( 0,  1,  4),  & v_8=( 0,  1,  -1),  & v_9=( 1,  2,  6),  \\
v_{10}=( 1,  2,  -1),  & v_{11}=( 1,  1,  -1),  & v_{12} =( 1,  3,  8),  \\
v_{13}=( 1,  3,  7),  & v_{14}=( 1,  3,  6),  & v_{15}=( 1,  3,  5),  \\
v_{16}=( 1,  3,  4),  & v_{17}=( 1,  3,  3),  & v_{18}=( 1,  3,  2),  \\
v_{19}=( 1,  3,  1),  & v_{20} =( 1,  3,  0),  
\end{array}
\]
and let $\tilde{D_i}$ be the toric divisor determined by the lattice point $v_i$ for $i=1,\ldots,20$, and $D_i:=\tilde{D_i}|_{{-}K_X}$ with $X:=\widetilde{\mathbb{P}_{\Sigma}}$. 
It can be easily seen by formulas (\ref{PicardNumber}) and (\ref{self-intersection}) that 
\[
\begin{array}{llll}
\rho_{\Delta} = 20-3 = 17, & 
D_1^2 = 2, \, D_2^2 = \cdots = D_4^2 = -2, \\
 & D_5^2 = 0, \, D_6^2 = \cdots = D_{20}^2 = -2. 
\end{array}
\]
Let $L$ be a lattice generated by divisors $\{D_i \}_{i=1}^{20}$. 
By solving the equation (\ref{ToricDivisorsRelations}), one sees that $\{$ $D_5$, $D_1$, $D_{11}$, $D_2$, $D_9$, $D_3$, $D_{12}$, $D_{13}$, $D_{14}$, $D_{15}$, $D_{16}$, $D_{17}$, $D_{18}$, $D_{19}$, $D_{20}$, $D_4$, $D_8$ $\}$ form a basis for $L$. 
By taking a new basis
\[
\left\{
\begin{array}{l}
D_5,\, D_5+D_{11},\, D_2-2D_5-D_{11},\, -D_1+D_2+2D_5+D_{11},\\
D_9,\, D_3,\, D_{12},\, D_{13},\, D_{14},\, D_{15},\, D_{16},\, D_{17},\,
D_{18},\, D_{19},\, D_{20},\, D_4,\, D_8 
\end{array}
\right\}, 
\]
one sees that the lattice $L$ is isometric to $U\oplus\tilde{L}$, where $\tilde{L}$ is of rank $15$ and of discriminant $-2$. 

Set one-simplices of $\Sigma'$ as follows: 
\[
\begin{array}{lll}
m_1=( -1,  0,  0), & m_2=( 2,  -1,  0), & m_3=( 0,  0,  1), \\
m_4=( -2,  4,  -1), & m_5=( -1, 3,  -1), & m_6=( -1,  2, 0), 
\end{array}
\]
and let $\tilde{D_i'}$ be the toric divisor determined by the lattice point $m_i$ for $i=1,\ldots,6$, and $D_i':=\tilde{D_i'}|_{{-}K_X}$ with $X:=\widetilde{\mathbb{P}_{\Sigma'}}$. 
It can be easily seen by formulas (\ref{PicardNumber}) and (\ref{self-intersection}) that 
\[
\begin{array}{llll}
\rho_{\Delta'} = 6-3 = 3, & 
D_1'^2 = 16, \, D_2'^2 = 6, \, D_3'^2 = 0, \\
&  D_4'^2 = D_5'^2 = D_6'^2 = -2. 
\end{array}
\]
Let $L'$ be a lattice generated by divisors $\{D_i' \}_{i=1}^{6}$. 
By solving the equation (\ref{ToricDivisorsRelations}), one sees that $\{$ $D_3'$, $D_6'$, $D_5'$ $\}$ form a basis for $L'$, with respect to which the intersection matrix of $L'$ is given by 
$\left(\begin{smallmatrix}
 0 & 1 & 0 \\
 1 & -2 & 0 \\
 0 & 0 & -2 \\
\end{smallmatrix}\right)$. 
By taking a new basis $\left\{D_3',\, D_3'+D_6',\, D_5' \right\}$, one sees that the lattice $L'$ is isometric to $U\oplus A_1$, which is a primitive sublattice of the $K3$ lattice. 
Therefore, $\Pic_{\!\Delta'}\simeq U\oplus A_1$. 

In cases 1 and 2, we obtain a lattice $L\simeq U\oplus\tilde{L}$, where $\tilde{L}$ is a lattice satisfying the assumption of Lemma~\ref{Lemma38-40}. 
Therefore, $L$ is a primitive sublattice of the $K3$ lattice, and that $\Pic_{\!\Delta}=L$ holds. 
Since $\discr{\Pic_{\!\Delta}}=\discr{(U\oplus\Pic_{\!\Delta'})}=2$, by Corollary~\ref{NikulinOrthogonal}, the relation 
$(\Pic_{\!\Delta})^\perp_{\Lambda_{K3}}\simeq U\oplus\Pic_{\!\Delta'}$ holds. 
Moreover, by Lemma~\ref{Nishiyama}, we have $\Pic_{\!\Delta}\simeq(U^{\oplus2}\oplus A_1)^\perp_{\Lambda_{K3}}\simeq U\oplus E_7\oplus E_8$. 

\subsection{Nos. 41, 42, and 43}
Set one-simplices of $\Sigma$ in terms of a basis of $M_{(1,1,3,4)}\otimes\mathbb{R}$ 
\[
(-1,1,0,0), \, (-3,0,1,0), \, (-4,0,0,1) :
\]
\[
\begin{array}{lll}
v_1=( -1,  2,  -1), & v_2=( -1,  -1,  -1), & v_3=( 8,  -1,  -1), \\
v_4=( 0,  -1,  1), & v_5=( -1,  -1,  1), & v_6=( 2,  1,  -1), \\
v_7=( 5,  0,  -1), & v_8=( 0,  -1,  -1), & v_9=( 1,  -1,  -1), \\
v_{10}=( 2,  -1,  -1), & v_{11}=( 3,  -1,  -1), & v_{12}=( 4,  -1,  -1), \\
v_{13}=( 5,  -1, -1), & v_{14}=( 6,  -1,  -1), & v_{15}=( 7,  -1,  -1), \\
v_{16}=( 4,  -1,  0), & v_{17}=( -1,  -1,  0), & v_{18}=( -1,  1,  -1), \\
v_{19}=( -1,  0,  -1), 
\end{array}
\]
and let $\tilde{D_i}$ be the toric divisor determined by the lattice point $v_i$ for $i=1,\ldots,19$, and $D_i:=\tilde{D_i}|_{{-}K_X}$ with $X:=\widetilde{\mathbb{P}_{\Sigma}}$. 
It can be easily seen by formulas (\ref{PicardNumber}) and (\ref{self-intersection}) that 
\[
\rho_{\Delta} = 19-3 = 16, \qquad
D_1^2 = 0, \quad D_2^2 = \cdots = D_{19} =-2. 
\]
Let $L$ be a lattice generated by divisors $\{D_i \}_{i=1}^{19}$. 
By solving the equation (\ref{ToricDivisorsRelations}), one sees that $\{D_6$, $D_1$, $D_4$, $D_{18}$, $D_{19}$, $D_2$, $D_8$, $D_9$, $D_{10}$, $D_{11}$, $D_{17}$, $D_5$, $D_{16}$, $D_3$, $D_{15}$, $D_{14}\}$ form a basis for $L$. 
By taking a new basis
\[
\left\{
\begin{array}{l}
D_1+D_4,\, D_1,\, D_1-D_4-D_5+D_6-D_{17}+D_{18},\, D_{17},\, D_{19},\\ D_2+D_{19},\, D_8,\, D_9,\, D_{10},\, D_{11},\\
3D_1+D_2-D_5+D_8+D_9+D_{10}+D_{11}+D_{18}+D_{19},\\
-3D_1+D_5-D_6,\, -D_1+D_{16},\, D_3,\, D_{15},\, D_{14}
\end{array}
\right\}, 
\]
one sees that the lattice $L$ is isometric to $U\oplus \tilde{L}$ with some lattice $\tilde{L}$. 
By a direct computation, one sees that $\sign{L}=(t_+,\,t_-)=(1,\, 15),\, \discr{L}=-3$, and $\rk{L}=16$, and thus, $\discr{\tilde{L}}=3$ and $\rk{\tilde{L}}=14$ hold. 
In particulat, the discriminant group $A_L$ of $L$ is isomorphic to $\mathbb{Z}\slash3\mathbb{Z}$, and $l(A_L)=1$. 
Therefore, one observes that 
\[
19-t_-=4\geq 0, \quad
3-t_+=2\geq 0,\quad
22-\rk{L}=6>1=l(A_{L})
\]
and by Corollary \ref{NikulinPrimitive}, $L$ is a primitive sublattice of the $K3$ lattice. 
Therefore, $\Pic_{\!\Delta}\simeq U\oplus\tilde{L}$ with $\discr{\tilde{L}}=3$ and $\rk{\tilde{L}}=14$. 

Set one-simplices of $\Sigma'$ in terms of a basis of $M_{(3,4,11,18)}\otimes\mathbb{R}$ 
\[
(-1,8,-1,-1), \, (0,-1,2,-1), \, (-1,-1,-1,1) :
\]
\[
\begin{array}{lll}
m_1=( 1,  0,  0), & m_2=( 0,  1,  0), & m_3=( 0,  0,  1), \\
m_4=( -1,  -3,  -4), & m_5=( 0,  -2,  -3), & m_6=( 0,  0,  -1), \\
m_7=( 0,  -1,  -2), 
\end{array}
\]
and let $\tilde{D_i'}$ be the toric divisor determined by the lattice point $m_i$ for $i=1,\ldots,7$, and $D_i':=\tilde{D_i'}|_{{-}K_X}$ with $X:=\widetilde{\mathbb{P}_{\Sigma'}}$. 
It can be easily seen by formulas (\ref{PicardNumber}) and (\ref{self-intersection}) that 
\[
\begin{array}{llll}
\rho_{\Delta'} = 7-3 = 4, & 
D_1'^2 = 0, \, D_2'^2 = 6, \, D_3'^2 = 12, \\
&  D_4'^2 = 0, \, D_5'^2 = D_6'^2 = D_7'^2 = -2. 
\end{array}
\]
Let $L'$ be a lattice generated by divisors $\{D_i' \}_{i=1}^{7}$. 
By solving the equation (\ref{ToricDivisorsRelations}), one sees that $\{$ $D_4'$, $D_5'$, $D_7'$, $D_6'$ $\}$ form a basis for $L'$. 
By taking a new basis $\left\{D_4',\, D_4'+D_5',\, -D_4'+D_7',\, D_6' \right\}$, one sees that the lattice $L'$ is isometric to $U\oplus A_2$, which is a primitive sublattice of the $K3$ lattice. 
Therefore, $\Pic_{\!\Delta'}\simeq U\oplus A_2$. 

Since $\discr{\Pic_{\!\Delta}}=\discr{U\oplus \tilde{L}}=-\discr{U^{\oplus 2}\oplus A_2}=3$, by Corollary~\ref{NikulinOrthogonal}, the relation $(\Pic_{\!\Delta})^\perp_{\Lambda_{K3}}\simeq U\oplus \Pic_{\!\Delta'}$ holds. 
Moreover, by Lemma~\ref{Nishiyama}, we have $\Pic_{\!\Delta}\simeq (U^{\oplus 2}\oplus A_2)^\perp_{\Lambda_{K3}}\simeq U\oplus E_6\oplus E_8$. 

\subsection{No. 46}
Set one-simplices of $\Sigma$ in terms of a basis of $M_{(1,1,1,2)}\otimes\mathbb{R}$  
\[
(-1,1,0,0),\, (-1,0,1,0),\, (-2,0,0,1) :
\]
\[
\begin{array}{lll}
v_1=(-1, -1, 1), & v_2=(-1, 2, 0), & v_3=(2, -1, 0), \\
v_4=(4,  -1,  -1),& v_5=(-1, 4, -1), & v_6=(-1, -1,  -1), \\
v_7=(-1, -1,  0), & v_8=(1, 0, 0), & v_9=(0, 1, 0), \\
v_{10}=(3, 0, -1), & v_{11}=(2, 1, -1), & v_{12}=(1, 2, -1), \\
v_{13}=(0, 3, -1), & v_{14}=(-1, 3, -1), & v_{15}=(-1, 2, -1), \\
v_{16}=(-1, 1, -1), & v_{17}=(-1, 0, -1), & v_{18}=(0, -1, -1), \\
v_{19}=(1, -1, -1), & v_{20}=(2, -1, -1), & v_{21}=(3, -1, -1), 
\end{array}
\]
and let $\tilde{D_i}$ be the toric divisor determined by the lattice point $v_i$ for $i=1,\ldots,21$, and $D_i:=\tilde{D_i}|_{{-}K_X}$ with $X:=\widetilde{\mathbb{P}_{\Sigma}}$. 
It can be easily seen by formulas (\ref{PicardNumber}) and (\ref{self-intersection}) that 
\[
\rho_{\Delta} = 21-3 = 18, \qquad
D_1^2 = 0, \quad D_2^2 =\cdots = D_{21}^2 = -2. 
\]
Let $L$ be a lattice generated by divisors $\{D_i \}_{i=1}^{21}$. 
By solving the equation (\ref{ToricDivisorsRelations}), one sees that $\{$ $D_1$,$D_3$,$D_8$,$D_7$,$D_2$,$D_5$,$D_{13}$,$D_{12}$,$D_{11}$,$D_{10}$,$D_{14}$,$D_{15}$,$D_{16}$,\\ 
$D_{17}$,$D_6$,$D_{18}$,$D_{19}$,$D_{20}$ $\}$ form a basis for $L$. 
Since $\rk{L}=18$ is strictly greater than $12$, the lattice $L$ is isometric to $U\oplus\tilde{L}$ with some lattice $\tilde{L}$. 
By a direct computation, one sees that $\sign{L}=(t_+,\,t_-)=(1,\, 17),\, \discr{L}=-5$, and $\rk{L}=18$, and thus, $\discr{\tilde{L}}=5$ and $\rk{\tilde{L}}=16$ hold. 
In particulat, the discriminant group $A_L$ of $L$ is isomorphic to $\mathbb{Z}\slash5\mathbb{Z}$, and $l(A_L)=1$. 
Therefore, one observes that 
\[
19-t_-=2\geq 0, \quad
3-t_+=2\geq 0,\quad
22-\rk{L}=4>1=l(A_{L})
\]
and by Corollary \ref{NikulinPrimitive}, $L$ is a primitive sublattice of the $K3$ lattice. 
Therefore, $\Pic_{\!\Delta}\simeq U\oplus\tilde{L}$ with $\discr{\tilde{L}}=5$ and $\rk{\tilde{L}}=16$. 

Set one-simplices of $\Sigma'$ in terms of a basis of $M_{(4,5,7,9)}\otimes\mathbb{R}$
\[
(4,0,-1,-1), \, (3,-1,-1,0), \, (0,-1,2,-1) :
\]
\[
\begin{array}{lll}
m_1=(0, 0, 1), & m_2=(2, -3, -1), & m_3=(-1, 1, 0), \\
m_4=(0, 1, 0), & m_5=(1, 0, 0), 
\end{array}
\]
and let $\tilde{D_i'}$ be the toric divisor determined by the lattice point $m_i$ for $i=1,\ldots,5$, and $D_i':=\tilde{D_i'}|_{{-}K_X}$ with $X:=\widetilde{\mathbb{P}_{\Sigma'}}$. 
It can be easily seen by formulas (\ref{PicardNumber}) and (\ref{self-intersection}) that 
\[
\rho_{\Delta'} = 5-3 = 2, \qquad
D_1'^2 = D_2'^2 = 2, \quad D_3'^2 = 10, \quad D_4'^2 = D_5'^2 = -2. 
\]
Let $L'$ be a lattice generated by divisors $\{D_i' \}_{i=1}^{5}$. 
By solving the equation (\ref{ToricDivisorsRelations}), one sees that $\{D_1', D_5'\}$ form a basis for $L$, with respect to which the intersection matrix of $L$ is given by $\left(\begin{smallmatrix} 2 & 1 \\ 1 & -2 \end{smallmatrix}\right)$. 
One sees that the lattice $L'$ is a hyperbolic lattice, that is, of signature $(t_+,\, t_-)=(1,1)$ of $\rk{L'}=2$ and $\discr{L'}=-5$. 
In particulat, the discriminant group $A_{L'}$ of $L'$ is isomorphic to $\mathbb{Z}\slash5\mathbb{Z}$, and $l(A_{L'})=1$. 
Therefore, one observes that 
\[
19-t_-=18\geq 0, \quad
3-t_+=2\geq 0,\quad
22-\rk{L'}=20>1=l(A_{L'})
\]
and by Corollary \ref{NikulinPrimitive}, $L'$ is a primitive sublattice of the $K3$ lattice. 
Therefore, $\Pic_{\!\Delta'}\simeq\left(\mathbb{Z}^2,\, \left(\begin{smallmatrix} 2 & 1 \\ 1 & -2 \end{smallmatrix}\right)\right)$. 

Since $\discr{\Pic_{\!\Delta}}=-\discr{\Pic_{\!\Delta}}=-5$, by Corollary~\ref{NikulinOrthogonal}, the relation $(\Pic_{\!\Delta})^\perp_{\Lambda_{K3}}\simeq U\oplus \Pic_{\!\Delta'}$ holds. 

\subsection{Nos. 48 and 49}
Set one-simplices of $\Sigma'$ in terms of a basis of $M_{(5,6,8,11)}\otimes\mathbb{R}$ 
\[
(-1,0,2,-1), \,  (-1,-1,0,1), \, (5,-1,-1,-1) :
\]
\[
\begin{array}{llll}
m_1=(1,0,0), & m_2=(0,1,0), & m_3=(0,0,1), & m_4=(-1,3,-1), 
\end{array}
\]
and let $\tilde{D_i'}$ be the toric divisor determined by the lattice point $m_i$ for $i=1,\ldots,4$, and $D_i':=\tilde{D_i'}|_{{-}K_X}$ with $X:=\widetilde{\mathbb{P}_{\Sigma'}}$. 
It can be easily seen by formulas (\ref{PicardNumber}) and (\ref{self-intersection}) that 
\[
\rho_{\Delta'} = 4-3 = 1, \qquad
D_1'^2 = D_2'^2 = 2, \quad D_3'^2 = 18, \quad D_4'^2 = 2. 
\]
Let $L'$ be a lattice generated by divisors $\{D_i' \}_{i=1}^{4}$. 
By solving the equation (\ref{ToricDivisorsRelations}), one sees that $\{D_1'\}$ form a basis for $L'$. 
Therefore, $\Pic_{\!\Delta'}\simeq\langle 2\rangle$. 
It is well-known that the lattice $\left( \mathbb{Z},\,\langle 2\rangle\right)$ is a primitive sublattice of the $K3$ lattice. 

Set one-simplices of $\Sigma$ in terms of a basis of $M_{(1,1,1,3)}\otimes\mathbb{R}$
\[
(-1,1,0,0), \, (-1,0,1,0), \, (-3,0,0,1) :
\]
\[
\begin{array}{lll}
v_1=( -1,  -1,  1), & v_2=( -1,  -1,  -1), & v_3=( 5, -1,  -1), \\
v_4=( -1,  -1,  5), & v_5=( -1,  -1,  0), & v_6=( 2,  -1,  0), \\
v_7=( -1,  -1,  2), & v_8=( 0,  -1,  -1), & v_9=( 1,  -1,  -1), \\
v_{10}=( 2,  -1,  -1), & v_{11}=( 3,  -1,  -1), & v_{12}=( 4,  -1,  -1), \\
v_{13}=( 4,  -1,  0), & v_{14}=( 3,  -1,  1), & v_{15}=( 2,  -1,  2), \\
v_{16}=( 1,  -1,  3), & v_{17}=( 0,  -1,  4), & v_{18}=( -1,  -1,  4), \\
v_{19}=( -1,  -1,  3), & v_{20}=( -1,  -1, 2), & v_{21}=( -1,  -1,  1), \\
v_{22}=( -1,  -1,  0), 
\end{array}
\]
and let $\tilde{D_i}$ be the toric divisor determined by the lattice point $v_i$ for $i=1,\ldots,22$, and $D_i:=\tilde{D_i}|_{{-}K_X}$ with $X:=\widetilde{\mathbb{P}_{\Sigma}}$. 
It can be easily seen by formulas (\ref{PicardNumber}) and (\ref{self-intersection}) that 
\[
\rho_{\Delta} = 22-3 = 19, \qquad
D_1^2 = 0, \quad D_2^2 = \cdots  = D_{22}^2 = -2. 
\]
Let $L$ be a lattice generated by divisors $\{D_i \}_{i=1}^{22}$. 
By solving the equation (\ref{ToricDivisorsRelations}), one sees that $\{$ $D_{21}$, $D_{20}$, $D_{19}$, $D_{18}$, $D_4$, $D_{17}$, $D_{16}$, $D_{15}$, $D_{14}$, $D_{13}$, $D_3$, $D_{12}$, $D_{11}$, $D_{10}$, $D_9$, $D_6$, $D_1$, $D_5$, $D_2$ $\}$ form a basis for $L$. 
By taking a new basis
\[
\left\{
\begin{array}{l}
D_1,\, D_1+D_5,\, D_{19},\, D_{18},\, D_4,\, D_{17},\, D_{16},\, D_{15},\, D_{14},\, D_{13},\, D_3,\, D_{12},\\
D_{11},\, D_{10},\, D_9,\, -D_1-D_2-D_5+D_6,\, D_{21},\, D_{20},\, D_2-D_1
\end{array}
\right\}, 
\]
one sees that the lattice $L$ is isometric to $U\oplus \tilde{L}$ with some lattice $\tilde{L}$. 

By a direct computation, one sees that $\sign{L}=(t_+,\,t_-)=(1,\, 18)$,\, $\discr{L}=2$, and $\rk{L}=19$, and thus, $\discr{\tilde{L}}=-2$ and $\rk{\tilde{L}}=17$ hold. 
In particulat, the discriminant group $A_L$ of $L$ is isomorphic to $\mathbb{Z}\slash2\mathbb{Z}$, and $l(A_L)=1$. 
Therefore, one observes that 
\[
19-t_-=1\geq 0, \quad
3-t_+=2\geq 0,\quad
22-\rk{L}=3>1=l(A_{L})
\]
and by Corollary \ref{NikulinPrimitive}, $L$ is a primitive sublattice of the $K3$ lattice. 
Therefore, $\Pic_{\!\Delta}\simeq U\oplus\tilde{L}$ with $\discr{\tilde{L}}=-2$ and $\rk{\tilde{L}}=17$. 

Since $\discr{\Pic_{\!\Delta}}=-\discr{(U\oplus\Pic_{\!\Delta'})}=2$, by Corollary~\ref{NikulinOrthogonal}, the relation $(\Pic_{\!\Delta})^\perp_{\Lambda_{K3}}\simeq U\oplus \Pic_{\!\Delta'}$ holds. 

Since the rank-one lattice $\left( \mathbb{Z},\, \langle 2\rangle\right)$ can be primitively embedded into the hyperbolic lattice $U=\langle e,\, f\rangle_{\mathbb{Z}}$ of rank $2$ as an element $e+f$, the orthogonal complement $\left( \mathbb{Z},\, \langle 2\rangle\right)^\perp_{U}$ in $U$ is a rank-one lattice $\left( \mathbb{Z},\, \langle -2\rangle\right)=\langle e-f \rangle_{\mathbb{Z}}$. 
Therefore, we have $\Pic_{\!\Delta'}\simeq \left( \mathbb{Z},\, \langle 2\rangle\right)^\perp_{U}\oplus U\oplus E_8^{\oplus 2}\simeq \left( \mathbb{Z},\, \langle -2\rangle\right)\oplus U\oplus E_8^{\oplus 2}$. 

\subsection{No. 50}
Set one-simplices of $\Sigma'$ in terms of a basis of $M_{(7,8,9,12)}\otimes\mathbb{R}$  
\[
 (-1,2,-1,0), \,  (-1,-1,3,-1), \,  (-1,-1,-1,2) :
\]
\[
\begin{array}{llll}
m_1=(1,0,0), & m_2=(0,1,0), & m_3=(0,0,1), & m_4=(-1,-1,-1), 
\end{array}
\]
and let $\tilde{D_i'}$ be the toric divisor determined by the lattice point $m_i$ for $i=1,\ldots,4$, and $D_i':=\tilde{D_i'}|_{{-}K_X}$ with $X:=\widetilde{\mathbb{P}_{\Sigma'}}$. 
One can easily seen by formulas (\ref{PicardNumber}) and (\ref{self-intersection}) that 
\[
\rho_{\Delta} = 4-3 = 1, \quad D_1'^2 = D_2'^2 = D_3'^2 = D_4'^2 = 4. 
\]
Let $L'$ be a lattice generated by divisors $\{D_i' \}_{i=1}^{4}$. 
By solving the equation (\ref{ToricDivisorsRelations}), one sees that $\{D_1'\}$ form a basis for $L'$. 
Therefore, $\Pic_{\!\Delta'}\simeq\left( \mathbb{Z},\, \langle 4\rangle\right)$. 
It is well-known that the lattice $\left( \mathbb{Z},\, \langle 4\rangle\right)$ is a primitive sublattice of the $K3$ lattice. 

Set one-simplices of $\Sigma$ in terms of a basis of $M_{(1,1,1,1)}\otimes\mathbb{R}$  
\[
 (-1,1,0,0), \,  (-1,0,1,0), \,  (-1,0,0,1) :
\]
\[
\begin{array}{lll}
v_1=( -1,  -1,  -1), & v_2=( 3,  -1,  -1), & v_3=( -1,  3,  -1), \\
v_4=( -1,  -1,  3), & v_5=( 0, -1,  -1), & v_6=( 1,  -1,  -1), \\
v_7=( 2,  -1,  -1), & v_8=( -1,  0,  -1), & v_9=( -1,  1,  -1), \\
v_{10}=( -1,  2,  -1), & v_{11}=( -1,  -1,  0), & v_{12}=( -1,  -1,  1), \\
v_{13}=( -1,  -1,  2), & v_{14}=( 2,  0,  -1), & v_{15}=( 1,  1,  -1), \\
v_{16}=( 0,  2,  -1), & v_{17}=( -1,  2,  0), & v_{18}=( -1,  1,  1), \\
v_{19}=( -1,  0,  2), & v_{20}=( 0,  -1,  2), & v_{21}=( 1,  -1,  1), \\
v_{22}=( 2,  -1,  0), 
\end{array}
\]
and let $\tilde{D_i}$ be the toric divisor determined by the lattice point $v_i$ for $i=1,\ldots,22$, and $D_i:=\tilde{D_i}|_{{-}K_X}$ with $X:=\widetilde{\mathbb{P}_{\Sigma}}$. 
It can be easily seen by formulas (\ref{PicardNumber}) and (\ref{self-intersection}) that 
\[
\rho_{\Delta} = 22-3 = 19, \quad D_1^2  = \cdots  = D_{22}^2 = -2. 
\]
Let $L$ be a lattice generated by divisors $\{D_i \}_{i=1}^{22}$. 
By solving the equation (\ref{ToricDivisorsRelations}), one sees that $\{$ $D_6$, $D_5$, $D_8$, $D_1$, $D_{11}$, $D_{12}$, $D_{13}$, $D_4$, $D_{19}$, $D_{18}$, $D_{17}$, $D_3$, $D_{10}$, $D_{16}$, $D_{15}$, $D_{14}$, $D_2$, $D_{22}$, $D_{21}$ $\}$ form a basis for $L$. 
Since the rank of $L$ is strictly greater than $12$, the lattice $L$ is isometric to $U\oplus \tilde{L}$ with some lattice $\tilde{L}$. 

By a direct computation, one sees that $\sign{L}=(t_+,\,t_-)=(1,\, 18),\, \discr{L}=4$, and $\rk{L}=19$, and thus, $\discr{\tilde{L}}=-4$ and $\rk{\tilde{L}}=17$ hold. 
In particulat, the discriminant group $A_L$ of $L$ is isomorphic to $\mathbb{Z}\slash4\mathbb{Z}$, and $l(A_L)=1$. 
Therefore, one observes that 
\[
19-t_-=1\geq 0, \quad
3-t_+=2\geq 0,\quad
22-\rk{L}=3>1=l(A_{L})
\]
and by Corollary \ref{NikulinPrimitive}, $L$ is a primitive sublattice of the $K3$ lattice. 
Therefore, $\Pic_{\!\Delta}\simeq U\oplus\tilde{L}$ with $\discr{\tilde{L}}=-4$ and $\rk{\tilde{L}}=17$. 

Since $\discr{\Pic_{\!\Delta}}=-\discr{(U\oplus\Pic_{\!\Delta'})}=2$, by Corollary~\ref{NikulinOrthogonal}, the relation $(\Pic_{\!\Delta})^\perp_{\Lambda_{K3}}\simeq U\oplus \Pic_{\!\Delta'}$ holds. 

Since the rank-one lattice $\left( \mathbb{Z},\, \langle 4\rangle\right)$ can be primitively embedded into the hyperbolic lattice $U=\langle e,\, f\rangle_{\mathbb{Z}}$ of rank $2$ as an element $2e+f$, the orthogonal complement $\left( \mathbb{Z},\, \langle 4\rangle\right)^\perp_{U}$ in $U$ is a rank-one lattice $\left( \mathbb{Z},\, \langle -4\rangle\right)=\langle e-2f \rangle_{\mathbb{Z}}$. 
Therefore, we have $\Pic_{\!\Delta'}\simeq \left( \mathbb{Z},\, \langle 4\rangle\right)^\perp_{U}\oplus U\oplus E_8^{\oplus 2}\simeq \left( \mathbb{Z},\, \langle -4\rangle\right)\oplus U\oplus E_8^{\oplus 2}$. 

Therefore, the assertion of Theorem~\ref{MainThm} is verified. 
\QED
\section{Conclusion}
We see in the main theorem that all coupling pairs that are polytope-dual with trivial toric contribution can extend to lattice duality among families of $K3$ surfaces. 
Thus, the coupling is partly translated to be the lattice-duality. 
Moreover, all except Nos. $46$, $48$ and $49$, and $50$ admit a pair of families of $K3$ surfaces with generic sections being elliptic: indeed, the Picard lattices $\Pic_{\!\Delta}$ and $\Pic_{\!\Delta'}$ contain the hyperbolic lattice $U$ of rank $2$. 

We can conclude that the Picard lattices of the families studied in the article are independent from the choice of reflexive polytopes. 
In other words, since the choice of a reflexive polytope is that of a way of blow-up of the ambient space, the Picard lattice in the subfamilies is birationally independent.

Makiko Mase\\
Universit\"at Mannheim, Lehrstuhl f\"ur Mathematik VI\\
B6, 26, 68131 Mannheim, Germany \\
email: mmase@mail.uni-mannheim.de
\end{document}